\documentclass[11pt]{article}%
\usepackage{amsmath}
\usepackage{graphicx}%
\usepackage{amsfonts}%
\usepackage{amssymb}
\newcommand{\eqnb}{\begin{equation}}
\newcommand{\eqne}{\end{equation}}

\newtheorem{The}{Theorem}

\newtheorem{Lem}{Lemma}

\newtheorem{Rem}{Remark}

\begin{document}

\title{Queueing Analysis of a Large-Scale Bike Sharing System through Mean-Field Theory}
\author{Quan-Lin Li, Chang Chen, Rui-Na Fan, Liang Xu and Jing-Yu Ma\\School of Economics and Management Sciences \\Yanshan University, Qinhuangdao 066004, P.R. China}
\maketitle

\begin{abstract}
The bike sharing systems are fast increasing as a public transport mode in
urban short trips, and have been developed in many major cities around the
world. A major challenge in the study of bike sharing systems is that some
large-scale and complex queueing networks have to be applied through
multi-dimensional Markov processes, while their discussion always suffers a
common difficulty: State space explosion. For this reason, this paper provides
a mean-field computational method to study such a large-scale bike sharing
system. Our mean-field computation is established in the following three
steps: Firstly, a multi-dimensional Markov process is set up for expressing
the states of the bike sharing system, and the empirical measure process of
the multi-dimensional Markov process is given to partly overcome the
difficulty of state space explosion. Based on this, the mean-field equations
are derived by means of a virtual time-inhomogeneous $M(t)/M(t)/1/K$ queue
whose arrival and service rates are determined by some mean-field computation.
Secondly, the martingale limit is employed to investigate the limiting
behavior of the empirical measure process, the fixed point is proved to be
unique so that it can be computed by means of a nonlinear birth-death process,
the asymptotic independence of this system is discussed, and specifically,
these lead to numerical computation for the steady-state probability of the
problematic (empty or full) stations. Finally, some numerical examples are
given for valuable observation on how the steady-state probability of the
problematic stations depends on some crucial parameters of the bike sharing system.

\vskip
0.5cm


\end{abstract}

\section{Introduction}

The bike sharing systems are fast developing wide-spread adoption in major
cities around the world, and are becoming a public mode of transportation
devoted to short trips. Up to now, there have been more than 500 cities
equipped with the bike sharing systems. Also, it is worth noting that the bike
sharing systems are being regarded as a promising solution to jointly
reducing, such as, traffic congestion, parking difficulty, transportation
noise, air pollution and global warming. For a history overview of the bike
sharing systems, readers may refer to, for instance, DeMaio \cite{DeM:2009}
and Shaheen et al. \cite{Sha:2010} for more details. While DeMaio
\cite{DeM:2003} provided a valuable prospect of the bike sharing systems in
the 21st century. For the status of the bike sharing systems in some countries
or cities, important examples include the United States by DeMaio and Gifford
\cite{DeM:2004}, France by Faye \cite{Fay:2008}, the European cities with OBIS
Project by Janett and Hendrik \cite{Jan:2011}, London by Lathia et al.
\cite{Lat:2012}, Montreal by Morency et al. \cite{Mor:2011}, Beijing by Liu et
al. \cite{Liu:2012}, several famous cities by Shu et al. \cite{Shu:2013}, and
further analysis by Shaheen et al. \cite{Sha:2011} and Meddin and DeMaio
\cite{Med:2012}.

To understand the recent key research directions, here it is necessary to
discuss some basic issues in design, operations and optimization of the bike
sharing systems. The literature of bike sharing systems may be classified into
two classes: The primary issues, and the higher issues. The primary issues are
to discuss the number of stations, the station location, the number of bikes,
the parking positions, and the types of bikes, all of which may be regarded as
the strategic design. The higher issues are to analyze the demand prediction,
the path scheduling, the inventory management, the repositioning (or
rebalancing) by trucks, the price incentive, and applications of the
intelligent information technologies. For analysis of the primary issues,
readers may refer to, for example, Dell'Olio et al. \cite{Del:2011}, Lin and
Yang \cite{Lin:2011}, Kumar and Bierlaire \cite{Kum:2012}, Martinez et al.
\cite{Mar:2012} and Nair et al. \cite{Nai:2013}. While the higher issues were
discussed by slightly more literature. Readers may refer to recent
publications or technical reports for more details, among which are
\textbf{the repositioning} by Forma et al. \cite{For:2010}, Vogel and Mattfeld
\cite{Vog:2010}, Benchimol et al. \cite{Ben:2011}, Raviv et al.
\cite{Rav:2013}, Contardo et al. \cite{Con:2012}, Caggiani and Ottomanelli
\cite{Cag:2012}, Fricker et al. \cite{Fri:2012}, Chemla et al. \cite{Che:2013}%
, Shu et al. \cite{Shu:2013}, Fricker and Gast \cite{Fri:2014} and Labadi et
al. \cite{Lab:2015}; \textbf{the inventory management} by Lin et al.
\cite{Lin:2013}, Raviv and Kolka \cite{Rav:2013} and Schuijbroek et al.
\cite{Sch:2013}; \textbf{the price incentives} by Waserhole and Jost
\cite{Was:2012}, Waserhole et al. \cite{Was:2013} and Fricker and Gast
\cite{Fri:2014}; \textbf{the fleet management} by George and Xia
\cite{Geo:2011, Geo:2010}, Godfrey and Powell \cite{God:2002}, Nair and
Miller-Hooks \cite{Nai:2011} and Guerriero et al. \cite{Gue:2012}; \textbf{the
simulation models} by Barth and Todd \cite{Bar:1999} and Fricker and Gast
\cite{Fri:2014}; \textbf{the data analysis} by Froehlich and Oliver
\cite{Fro:2008}, Vogel et al \cite{Vog:2011}, Borgnat et al. \cite{Bor:2011},
C\^{o}me et al. \cite{Com:2012} and Katzev \cite{Kat:2003}.

Based on the above literature, it is necessary to further observe a basic
solution to operations of the bike sharing systems. In a bike sharing system,
a customer arrives at a station, takes a bike, and uses it for a while; then
he returns the bike to a destination station. In general, the bikes are
frequently distributed in an imbalanced manner among the stations, thus an
arriving customer may always be confronted with two problematic cases: (1) A
station is empty when a customer arrives at the station to rent a bike, and
(2) a station is full when a bike-riding customer arrives at the station to
return his bike. For the two problematic cases, the empty or full station is
called a problematic station. Since a crucial question for the operational
efficiency of the bike sharing system is its ability not only to meet the
fluctuating demand for renting bikes at each station but also to provide
enough vacant lockers to allow the renters to return the bikes at their
destinations, the two types of problematic stations reflect a common challenge
facing operations management of the bike sharing systems in practice due to
the stochastic and time-inhomogeneous nature of the customer arrivals and bike
returns. Therefore, it is a key to measuring the steady-state probability of
the problematic stations in the study of bike sharing systems. Also, analysis
of the steady-state probability of the problematic stations is useful and
helpful in design, operations and optimization of the bike sharing systems in
terms of numerical computation and comparison. Up to now, it is still
difficult (and even impossible) to provide an explicit expression for the
steady-state probability of the problematic stations because the bike sharing
system is a more complicated closed queueing network with various geographical
interactions, which come both from some bikes parked in multiple stations and
from the other bikes ridden on multiple roads. For this, Section 2 explains
that the bike sharing system is a Markov process of dimension $N^{2}$ through
analysis of a complicated virtual closed queueing network, also see Li et al.
\cite{LiF:2016} for more details.

To compute the steady-state probability of the problematic stations, it is
better to develop a stochastic and dynamic method through applications of the
queueing theory as well as Markov processes to the study of bike sharing
systems. However, the available works on such a research direction are still
few up to now. To survey the recent literature, some significant methods and
results are listed as follows. \textbf{The simple queues: }Leurent
\cite{Leu:2012} used the $M/M/1/C$ queue to consider a vehicle-sharing system
in which each station contains an expanded waiting room only for those
customers arriving at either a full station to return a bike or an empty
station to rent a bike, and analyzed performance measures of this
vehicle-sharing system in terms of a geometric distribution. Schuijbroek et
al. \cite{Sch:2013} first computed the transient distribution of the $M/M/1/C$
queue, which is used to measure the service level in order to establish a
mixed integer programming for the bike sharing system. Then they dealt with
the inventory rebalancing and the vehicle routing by means of the optimal
solution to the mixed integer programming. Raviv et al \cite{Rav:2013} and
Raviv and Kolka \cite{Rav:2013a} provided an effective method for computing
the transient distribution of a time-inhomogeneous $M\left(  t\right)
/M\left(  t\right)  /1/C$ queue, which is used to evaluate the expected number
of bike shortages at any station. \textbf{The queueing networks: }Savin et al.
\cite{Sav:2005} used a loss network as well as the admission control to
discuss capacity allocation of a rental model with two classes of customers,
and studied the revenue management and fleet sizing decision. Adelman
\cite{Ade:2007} applied a closed queueing network to propose an internal
pricing mechanism for managing a fleet of service units, and also used a
nonlinear flow model to discuss the price-based policy for the vehicle
redistribution. George and Xia \cite{Geo:2011} provided an effective method of
closed queueing networks in the study of vehicle rental systems, and
determined the optimal number of parking spaces for each rental location. Li
et al. \cite{LiF:2016} proposed a unified framework for analyzing the closed
queueing networks in the study of bike sharing systems. \textbf{The mean-field
theory: }Recently, the mean-field method as well as the queueing theory are
applied to analyzing the bike sharing systems. Fricker et al. \cite{Fri:2012}
considered a space-inhomogeneous bike sharing system, and expressed the
minimal proportion of problematic stations within each cluster. Fricker and
Gast \cite{Fri:2014} provided a detailed analysis for a space-homogeneous bike
sharing system in terms of the $M/M/1/K$ queue and some simple mean-field
models, and crucially, they gave the closed-form solution to the minimal
proportion of problematic stations. Fricker and Tibi \cite{Fri:2015} first
studied the central limit and local limit theorems for the independent (non
identically distributed) random variables, which support analysis of a
generalized Jackson network with product-form solution; then they used the
limit theorems to give a better outline of the stationary asymptotic analysis
of the locally space-homogeneous bike sharing systems. Li and Fan
\cite{LiFe:2016} developed numerical computation of the bike sharing systems
under Markovian environment by means of the mean-field theory and the
nonlinear QBD processes. \textbf{The Markov decision processes: }A simple
closed queuing network is used to establish the Markov decision model in the
study of bike sharing systems, and to provide a fluid approximation in order
to compute the static optimal policy. Examples include Waserhole and Jost
\cite{Was:2012}, Waserhole and Jost \cite{Was:2013, Was:2014} and Waserhole et
al. \cite{Was:2013a}.\ \

For convenience of readers, it is necessary to recall some basic references in
which the mean-field theory is applied to the analysis of large-scale
stochastic systems. Readers may refer to Spitzer \cite{Spi:1970}, Dawson
\cite{Daw:1983}, Sznitman \cite{Szn:1989}, Vvedenskaya et al. \cite{Vve:1996},
Mitzenmacher \cite{Mit:1996}, Turner \cite{Tur:1998}, Graham \cite{Gra:2000,
Gra:2004}, Benaim and Le Boudec \cite{Ben:2008}, Gast and Gaujal
\cite{Gas:2010, Gas:2011}, Bordenave et al. \cite{Bor:2010}, Li \cite{Li:2014,
Li:2015a}, Li and Lui \cite{Li:2014b}, Li et al. \cite{Li:2014a, Li:2015},
Fricker et al. \cite{Fri:2012} and Fricker and Tibi \cite{Fri:2015}. On the
other hand, the metastability of Markov processes may be useful in the study
of more general bike sharing systems when the nonlinear Markov processes are
applied. Readers may refer to, such as, Bovier \cite{Bov:2003}, Den Hollander
\cite{Den:2004}, Antunes et al. \cite{Ant:2008}, Tibi \cite{Tib:2011}, Li
\cite{Li:2015a} and more references therein.

The main contributions of this paper are twofold. The first contribution is to
describe a mean-field queueing model to analyze the large-scale bike sharing
systems, where the arrival, walk, bike-riding (or return) processes among the
stations are given some simplified assumptions whose purpose is to guarantee
applicability of the mean-field theory. For this, we develop a mean-field
queueing method combining the mean-field theory with the time-inhomogeneous
queue, the martingale limits and the nonlinear birth-death processes. To this
end, we provide a complete picture of applying the mean-field theory to the
study of bike sharing systems through four basic steps: (1) \textbf{The system
of mean-field equations} is set up by means of a virtual time-inhomogeneous
$M(t)/M(t)/1/K$ queue whose arrival and service rates are determined by means
of some mean-field computation; (2) \textbf{the asymptotic independence} (or
\textbf{propagation of chaos}) is proved in terms of the martingale limit and
the uniqueness of the fixed point; (3) \textbf{numerical computation of the
fixed point} is given by using a system of nonlinear equations corresponding
to the nonlinear birth-death processes; and (4) \textbf{performance analysis}
of the bike sharing system is given through some numerical computation.

The second contribution of this paper is to provide a detailed analysis for
computing the steady-state probability of the problematic stations, which is
one of the most key measures in the study of bike sharing systems. It is worth
noting that the service level, optimal design and control mechanism of bike
sharing systems can be computed by means of the steady-state probability of
the problematic stations. Therefore, this paper develops effective algorithms
for computing the steady-state probability of the problematic stations, and
gives a numerically computational framework in the study of bike sharing
systems. Furthermore, we use some numerical examples to give valuable
observation and understanding on how the performance measures depend on some
crucial parameters of the bike sharing system. On the other hand, in view that
Fricker et al. \cite{Fri:2012}, Fricker and Gast \cite{Fri:2014}, Fricker and
Tibi \cite{Fri:2015} and Li and Fan \cite{LiFe:2016} are the only important
references that are closely related to this paper by using the mean-field
theory, but differently, this paper provides more work focusing on some key
theoretical points such as the virtual time-inhomogeneous $M(t)/M(t)/1/K$
queue, the mean-field equations, the martingale limits, the nonlinear
birth-death processes, numerical computation of the fixed point, and numerical
analysis for the steady-state probability of the problematic stations. With
successful exposition of the key theoretical points, such a numerical
computation can greatly enable a broad study of bike sharing systems.
Therefore, the methodology and results of this paper gain new insights on how
to establish the mean-field queueing models for discussing more general bike
sharing systems by means of the mean-field theory, the time-inhomogeneous
queues and the nonlinear Markov processes.

The remainder of this paper is organized as follows. In Section 2, we first
describe a large-scale bike sharing system with $N$ identical stations, give a
$N$-dimensional Markov process for expressing the states of the bike sharing
system, and establish an empirical measure process of the $N$-dimensional
Markov process in order to partly overcome the difficulty of state space
explosion. In Section 3, we set up a system of mean-field equations satisfied
by the expected fraction vector through a virtual time-inhomogeneous
$M(t)/M(t)/1/K$ queue whose arrival and service rates are determined by means
of some mean-field computation. In Section 4, we establish a Lipschitz
condition, and prove the existence and uniqueness of solution to the system of
mean-field equations. In Section 5, we provide a martingale limit of the
sequences of empirical measure Markov processes in the bike sharing system. In
Section 6, we analyze the fixed point of the system of mean-field equations,
and prove that the fixed point is unique. Based on this, we simply analyze the
asymptotic independence of the bike sharing system, and also discuss the
limiting interchangeability with respect to $N\rightarrow\infty$ and
$t\rightarrow+\infty$. In Section 7, we provide some effective computation of
the fixed point, and use some numerical examples to investigate how the
steady-state probability of the problematic stations depends on some crucial
parameters of the bike sharing system. Some concluding remarks are given in
Section 8.

\section{Model Description}

In this section, we first describe a large-scale bike sharing system with $N$
identical stations, and establish a $N$-dimensional Markov process for
expressing the states of the bike sharing system. To overcome the difficulty
of state space explosion, we provide an empirical measure process of the
$N$-dimensional Markov process.

We first show that a bike sharing system can be modeled as a complex
stochastic system whose analysis is always difficult and challenging. Then we
explain the reasons why it is necessary to develop some simplified models in
the study of bike sharing systems. In particular, we indicate that the
mean-field theory plays a key role in establishing and analyzing such a
simplified model whose purpose is to be able to set up some basic and useful
relations among several key parameters of system.

\textbf{A Complex Stochastic System}

In the bike sharing system, a customer arrives at a station, takes a bike, and
uses it for a while; then he returns the bike to any station and immediately
leaves this system. Based on this, if the bike sharing system has $N$ stations
for $N\geq2$, then it can contain at most $N\left(  N-1\right)  $ roads
because there may be a road between any two stations. When the stations and
roads are different and heterogeneous, Li et al. \cite{LiF:2016} showed that
the bike sharing system can be modeled as a complicated closed queueing
network due to the fact that the total number of bikes is fixed in this
system. In this case, the bikes are regarded as the virtual customers, while
the stations and the roads are viewed as the virtual servers. Based on this,
the closed queueing network is described as a Markov process $\left\{
\overrightarrow{\mathfrak{n}}\left(  t\right)  :t\geq0\right\}  $ of dimension
$N^{2}$, where%
\[
\overrightarrow{\mathfrak{n}}\left(  t\right)  =\left(  \mathbf{n}_{1}\left(
t\right)  ,\mathbf{n}_{2}\left(  t\right)  ,\ldots,\mathbf{n}_{N}\left(
t\right)  \right)  ,
\]%
\[
\mathbf{n}_{k}\left(  t\right)  =\left(  n_{k}\left(  t\right)  ;n_{k,1}%
\left(  t\right)  ,\ldots,n_{k,k-1}\left(  t\right)  ,n_{k,k+1}\left(
t\right)  ,\ldots,n_{k,N}\left(  t\right)  \right)  ,
\]%
\[
\sum_{k=1}^{N}n_{k}\left(  t\right)  +\sum_{i=1}^{N}\sum_{j\neq i}^{N}%
n_{i,j}\left(  t\right)  =\Re,
\]
$n_{k}\left(  t\right)  $ is the number of bikes parked at Station $k$,
$n_{k,j}\left(  t\right)  $ is the number of bikes rided on Road $k\rightarrow
j$ for $j\neq k$ and $1\leq j,k\leq N$, and $\Re$ is the total number of bikes
in the bike sharing system.

In general, analysis of the Markov process $\left\{  \overrightarrow
{\mathfrak{n}}\left(  t\right)  :t\geq0\right\}  $ of dimension $N^{2}$ is
usually difficult due to at least three reasons: (1) The state space explosion
for a large integer $N$, (2) the complex routes among the virtual servers
which are either the $N$ stations or the $N\left(  N-1\right)  $ roads, and
(3) a complicated expression for the steady-state probability distribution of
joint queue lengths. See Li and Fan \cite{LiF:2016} for more details. For
this, it is necessary in practice to provide a simplified model that contains
only several key parameters of system, while the simplified model is used to
set up some basic and useful relations among the key parameters. Crucially,
not only do the basic relations support numerical computation of the
steady-state probability of the problematic stations, but they are also
helpful for performance analysis of the bike sharing system. To provide such a
simplified model, the remainder of this paper will provide a mean-field
queueing model described from the bike sharing system.

\textbf{A Basic Condition to Apply the Mean-Field Theory}

To apply the mean-field theory, we only need to consider the bike information
$(n_{1}\left(  t\right)  ,n_{2}\left(  t\right)  $, $\ldots,n_{N}\left(
t\right)  )$ on the $N$ stations, while the bike information of the $N\left(
N-1\right)  $ roads will be combined into the `\textit{probabilistic
behavior}' of the random vector $\left(  n_{1}\left(  t\right)  ,n_{2}\left(
t\right)  ,\ldots,n_{N}\left(  t\right)  \right)  $ by means of some
mean-field computation. See Theorem \ref{The:Rates} and its proof in the next
section. At the same time, a basic condition is also needed to guarantee the
exchangeability of the $N$-dimensional Markov process $\left\{  \left(
n_{1}\left(  t\right)  ,n_{2}\left(  t\right)  ,\ldots,n_{N}\left(  t\right)
\right)  :t\geq0\right\}  $, that is, for any permutation $\left(  i_{1}%
,i_{2},i_{3},\ldots,i_{N}\right)  $ of $\left(  1,2,3,\ldots,N\right)  $,
\[
P\left\{  n_{1}\left(  t\right)  =k_{1},n_{2}\left(  t\right)  =k_{2}%
,\ldots,n_{N}\left(  t\right)  =k_{N}\right\}  =P\left\{  n_{i_{1}}\left(
t\right)  =k_{i_{1}},n_{i_{2}}\left(  t\right)  =k_{i_{2}},\ldots,n_{i_{N}%
}\left(  t\right)  =k_{i_{N}}\right\}  .
\]
See Li \cite{Li:2015a} for the mean-field analysis of big networks. In fact,
the following assumption \textbf{(1)} that the bike sharing system consists of
$N$ identical stations guarantee the exchangeability of the Markov process
$\left\{  \left(  n_{1}\left(  t\right)  ,n_{2}\left(  t\right)  ,\ldots
,n_{N}\left(  t\right)  \right)  :t\geq0\right\}  $ so that the mean-field
theory can be applied to discussing the bike sharing system.

Although the model assumptions to apply mean-field theory are simplified
greatly, we can still set up some useful and basic relations among several key
parameters\ of system, and also provide some simple and effective algorithms
both for computing the steady-state probability of the problematic stations
and for analyzing performance measures of the bike sharing system.

\textbf{Simplified Model Assumptions}

Based on the above analysis, we make some necessarily simplified assumptions
for applying the mean-field theory to studying the bike sharing system as follows:

\textbf{(1) The }$N$\textbf{ identical stations:} The bike sharing system
consists of $N$ identical stations, each of which has a finite bike capacity.
At the initial time $t=0$, each station contains $C$ bikes and $K$ positions
to park the bikes, where $1\leq C<K<\infty$.

\textbf{(2) The arrive processes: }The arrivals of outside customers at the
bike sharing system are a Poisson process with arrival rate $N\lambda$ for
$\lambda>0$.

\textbf{(3) The walk processes:} If an outside or walking customer arrives at
an empty station in which no bike may be rented, then he has to walk to
another station again in the hope of renting a bike. We assume that the
customer may rent a bike from a station within at most $\omega$ consequent
walks, otherwise he will directly leave this system (that is, if he has not
rented a bike after $\omega$ consequent walks yet). Note that one walk is
viewed as a process that the customer walks from an empty station to another
station, and $\omega$ is the maximal number of consequent walks of the
customer among the stations.

We assume that the walk times between any two stations are all exponential
with walk rate $\gamma>0$. Obviously, the expected walk time is $1/\gamma$.

\textbf{(4) The bike-riding (or return) processes:} If a bike-riding customer
arrives at a full station in which no parking position is available, then he
has to ride the bike to another station again. We assume that the
returning-bike process is persistent in the sense that the customer must find
a station with an empty position to return his bike (that is, he can not leave
this system before his bike is returned), because the bike is the public
property so that no one can make it his own.

We assume that the bike-riding times between any two stations are all
exponential with bike-riding rate $\mu$ for $\gamma\leq\mu<+\infty$. Clearly,
the expected bike-riding time is $1/\mu$.

\textbf{(5) The departure discipline: }The customer departure has two
different cases: (a)\textbf{ }The customer directly leaves the bike sharing
system if he has not rented a bike yet after $\omega$ consequent walks; or (b)
once one customer takes, uses and returns the bike to a station, he completes
this trip, thus he can immediately leave the bike sharing system.

We assume that the arrival, walk and bike-riding processes are independent,
and all the above random variables are independent of each other. Note that
the randomly bike-riding and walk times show that the road length between any
two stations is considered in this paper. For such a bike sharing system,
Figure 1 provides some physical interpretation.

\begin{figure}[ptb]
\centering             \includegraphics[width=8cm]{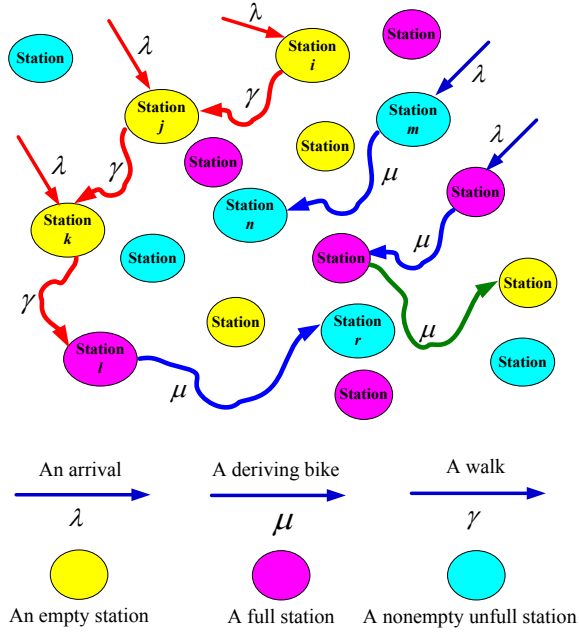}  \newline
\caption{The physical interpretation of a bike sharing system}%
\label{figure:figure-1}%
\end{figure}

\begin{Rem}
(1) The assumption of the $N$ identical stations is used to guarantee
applicability of the mean-field theory, that is, the $N$-dimensional Markov
process $\{(n_{1}\left(  t\right)  ,n_{2}\left(  t\right)  ,\ldots$,
$n_{N}\left(  t\right)  ):t\geq0\}$ is exchangeable. Although the model
assumptions to apply the mean-field theory are simplified greatly (note that
several key parameters of system will be observed and analyzed in a simple
form), we can still set up some useful and basic relations among the key
parameters\ of system, and also find some valuable law and pattern both from
computing the steady-state probability of the problematic stations and from
analyzing performance measures of the bike sharing system.

(2) It is necessary to explain the maximal number $\omega$ of consequent walks
of the customer. If $\omega=0$, then the arriving customer immediately leaves
this system once he arrives at a full station. If $\omega$ is smaller, then
the customer would like to find an available bike at a lucky station through
at most $\omega$ consequent walks, because a bike can help him to promptly
deal with a number of important things so that he would like to accept the
time delay due to the hope of renting a bike within at most $\omega$
consequent walks.

(3) The road lengths among the $N$ stations are considered here, while the
bike-riding time on any road is exponential with bike-riding rate $\mu$. Based
on this, the road length is measured by means of the randomly bike-riding
time. In addition, the assumption with $0<\gamma<\mu<+\infty$ makes sense\ in
practice because the riding bike is faster than the walk on any road. On the
other hand, the assumptions on the i.i.d. bike-riding times and on the i.i.d.
walk times are to guarantee applicability of the mean-field theory, that is,
the $N$-dimensional Markov process $\left\{  \left(  n_{1}\left(  t\right)
,n_{2}\left(  t\right)  ,\ldots,n_{N}\left(  t\right)  \right)  :t\geq
0\right\}  $ is exchangeable. Therefore, the $N$ identical stations also
contain their identically physical factors under a random setting.
\end{Rem}

In the remainder of this section, we first establish a $N$-dimensional Markov
process for expressing the states of the bike sharing system. Then we give an
empirical measure process of the $N$-dimensional Markov process in order to
overcome the difficulty of state space explosion.

Let $X_{i}^{\left(  N\right)  }\left(  t\right)  $ be the number of bikes
parked in Station $i$ at time $t\geq0$. Then $X_{i}^{\left(  N\right)
}\left(  t\right)  =n_{i}\left(  t\right)  $, and henceforth we only use the
notation $X_{i}^{\left(  N\right)  }\left(  t\right)  $. It is easy to see
from the above model descriptions that $\mathcal{X}=\left\{  \left(
X_{1}^{\left(  N\right)  }\left(  t\right)  ,X_{2}^{\left(  N\right)  }\left(
t\right)  ,\ldots,X_{N}^{\left(  N\right)  }\left(  t\right)  \right)
:t\geq0\right\}  $ is a $N$-dimensional Markov process. In general, it is
always more difficult to directly study the $N$-dimensional Markov process
$\mathcal{X}$ due to the state space explosion. Thus we need to introduce an
empirical measure process of the $N$-dimensional Markov process $\mathcal{X}$
as follows. We write%
\[
Y_{k}^{\left(  N\right)  }\left(  t\right)  =\frac{1}{N}\sum_{i=1}%
^{N}\mathbf{1}_{\left\{  X_{i}^{\left(  N\right)  }\left(  t\right)
=k\right\}  },
\]
where $\mathbf{1}_{\left\{  \cdot\right\}  }$ is an indicator function.
Obviously, $Y_{k}^{\left(  N\right)  }\left(  t\right)  $ is the proportion of
the stations with $k$ bikes at time $t$, and $0\leq\sum_{i=1}^{N}%
\mathbf{1}_{\left\{  X_{i}^{\left(  N\right)  }\left(  t\right)  =k\right\}
}\leq N$. Let%
\[
\mathbf{Y}^{\left(  N\right)  }\left(  t\right)  =\left(  Y_{0}^{\left(
N\right)  }\left(  t\right)  ,Y_{1}^{\left(  N\right)  }\left(  t\right)
,...,Y_{K-1}^{\left(  N\right)  }\left(  t\right)  ,Y_{K}^{\left(  N\right)
}\left(  t\right)  \right)  .
\]
Then it is easy to see that the empirical measure process $\left\{
\mathbf{Y}^{\left(  N\right)  }\left(  t\right)  :t\geq0\right\}  $ is a
Markov process on the state space $\Omega=\left[  0,1\right]  ^{K+1}$.

To study the empirical measure Markov process, we write%
\[
y_{k}^{\left(  N\right)  }\left(  t\right)  =E\left[  Y_{k}^{\left(  N\right)
}\left(  t\right)  \right]
\]
and%
\[
\mathbf{y}^{\left(  N\right)  }\left(  t\right)  =\left(  y_{0}^{\left(
N\right)  }\left(  t\right)  ,y_{1}^{\left(  N\right)  }\left(  t\right)
,...,y_{K-1}^{\left(  N\right)  }\left(  t\right)  ,y_{K}^{\left(  N\right)
}\left(  t\right)  \right)  .
\]

\section{The Mean-Field Equations}

In this section, we first describe the bike sharing system as a virtual
time-inhomogeneous $M(t)/M(t)/1/K$ queue whose arrival and service rates are
determined by means of the mean-field theory. Then we set up a system of
mean-field equations, which is satisfied by the expected fraction vector
$\mathbf{y}^{\left(  N\right)  }\left(  t\right)  $, in terms of the virtual
time-inhomogeneous $M(t)/M(t)/1/K$ queue.

Note that the $N$ stations are identical according to the above model
description on both system parameters and operations discipline, thus we can
use the mean-field theory to study the bike sharing system. In this case, we
only need to observe a tagged station (for example, Station 1) whose number of
bikes is regarded as a virtual time-inhomogeneous $M(t)/M(t)/1/K$ queue (see
Figure 2); while the other $N-1$ stations have some impact on the tagged
station, and the impact can be analyzed by means of the empirical measure
process through a mean-field computation for the new arrival and service rates
in this virtual queue. Specifically, we also explain the reason why the new
arrival and service processes in this virtual queue are time-inhomogeneous.
See Figure 2 for more details.

\begin{figure}[ptb]
\centering                 \includegraphics[width=10cm]{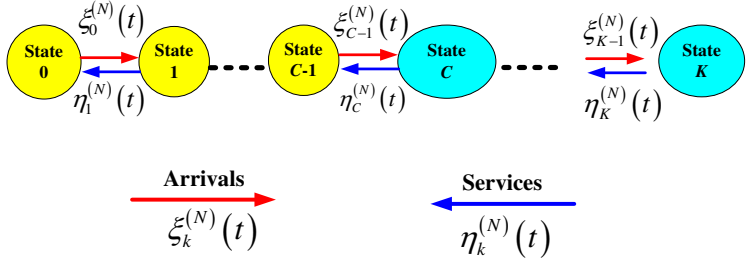}
\newline \caption{The state transitions in the $M(t)/M(t)/1/K$ queue}%
\label{figure:figure-2}%
\end{figure}

It is necessary to explain the difference of the arrival and service processes
between the bike sharing system and the virtual time-inhomogeneous
$M(t)/M(t)/1/K$ queue. For example, if a real customer arrives and rents a
bike at a tagged station, then the number of bikes parked in the tagged
station decreases by one, thus the real customer arrivals at the tagged
station should be a part of the service process of the $M(t)/M(t)/1/K$ queue;
while if a real customer returns a bike to a tagged station and leaves this
system (i.e., his trip is completed), then the number of bikes parked in the
tagged station increases by one, thus the real customers' returning their
bikes to the tagged station should be a part of the arrival process of the
$M(t)/M(t)/1/K$ queue. Furthermore, the following Theorem 1 provides a more
detailed analysis for various parts of the arrival and service processes in
the virtual time-inhomogeneous $M(t)/M(t)/1/K$ queue.

For the time-inhomogeneous $M(t)/M(t)/1/K$ queue, now we use the mean-field
theory to discuss its Poisson input with arrival rate $\xi_{l}^{\left(
N\right)  }\left(  t\right)  $ for $0\leq l\leq K-1$ and its exponential
service times with service rate $\eta_{k}^{\left(  N\right)  }\left(
t\right)  $ for $1\leq k\leq K$.

The following theorem provides expressions for the arrival and service rates:
$\xi_{l}^{\left(  N\right)  }\left(  t\right)  $ for $0\leq l\leq K-1$ and
$\eta_{k}^{\left(  N\right)  }\left(  t\right)  $ for $1\leq k\leq K$,
respectively. Note that the time-inhomogeneous arrival and service rates will
play a key role in our mean-field study later.

\begin{The}
\label{The:Rates}For $1\leq k\leq K$ and $\omega=0,1,2,\ldots$, we have%
\begin{equation}
\eta_{k}^{\left(  N\right)  }\left(  t\right)  =\eta^{\left(  N\right)
}\left(  t\right)  =\lambda+\gamma y_{0}^{\left(  N\right)  }\left(  t\right)
\frac{1-\left[  y_{0}^{\left(  N\right)  }\left(  t\right)  \right]  ^{\omega
}}{1-y_{0}^{\left(  N\right)  }\left(  t\right)  }. \label{Equ-1}%
\end{equation}
At the same time, for $0\leq l\leq K-1$ we have%
\begin{equation}
\xi_{l}^{\left(  N\right)  }\left(  t\right)  =\left\{
\begin{array}
[c]{ll}%
\frac{\mu}{N}\frac{1}{1-y_{K}^{\left(  N\right)  }\left(  t\right)  }\left\{
\left(  C-l\right)  +\left(  N-1\right)  \left[  C-\sum_{k=1}^{K}%
ky_{k}^{\left(  N\right)  }\left(  t\right)  \right]  \right\}  , & 0\leq
l\leq C-1,\\
\frac{\mu}{N}\frac{1}{1-y_{K}^{\left(  N\right)  }\left(  t\right)  }\left\{
\left(  N-1\right)  \left[  C-\sum_{k=1}^{K}ky_{k}^{\left(  N\right)  }\left(
t\right)  \right]  \right\}  , & C\leq l\leq K-1.
\end{array}
\right.  \label{Equ-2}%
\end{equation}
\end{The}

\textbf{Proof} \ We first prove Equation (\ref{Equ-1}). In this case, we need
to specifically deal with State $0$. If one customer arrives at an empty
station, then the customer has to walk from the empty station to another
station. It is easy to see that the bikes parked at the tagged station will
have two different cases: (a) There is at least one bike with probability
$\sum_{i=1}^{K}y_{i}^{\left(  N\right)  }\left(  t\right)  $; and (b) there is
no bike with probability $y_{0}^{\left(  N\right)  }\left(  t\right)  $. For
Case (a), the customer can rent a bike for his trip; while for Case (b), the
customer will have to walk to another station again until he hopes to be able
to rent a bike from a next station within the $\omega$ consequent walks.
Notice that the role played by State $0$ is depicted in Figure 3, thus we can
easily observe that the state transitions from State $0$ are jointly caused by
the arrival, walk and return (or bike-riding) processes.

\begin{figure}[ptb]
\centering   \includegraphics[width=9cm]{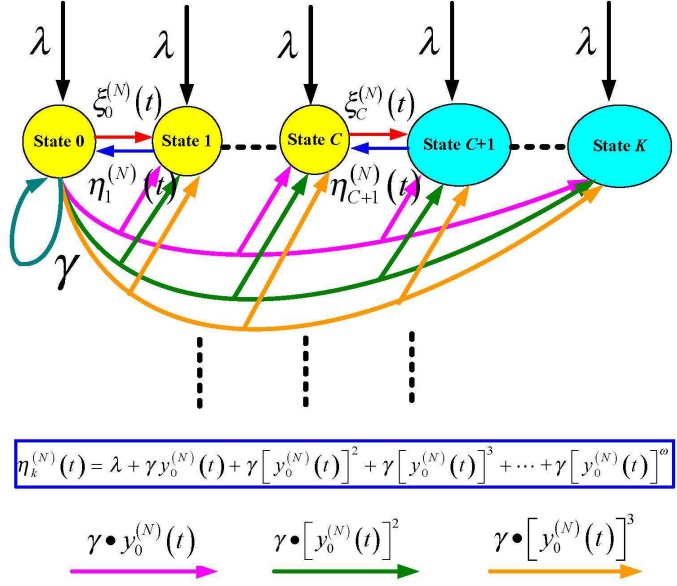}  \newline \caption{The
state transitions for computing $\eta_{k}^{\left(  N\right)  }\left(
t\right)  $}%
\label{figure:figure-3}%
\end{figure}

To compute the service rate $\eta_{k}^{\left(  N\right)  }\left(  t\right)  $
for $1\leq k\leq K$, it is seen from Figure 3 that State $0$ (that is, the
tagged station is empty) is a key, and it leads to the rate $\gamma\left[
y_{0}^{\left(  N\right)  }\left(  t\right)  \right]  ^{n}$ with respect to $n$
consequent walks, where the $n$ consequent\ walks correspond to $n$ empty
stations with probability $\left[  y_{0}^{\left(  N\right)  }\left(  t\right)
\right]  ^{n}$for $1\leq n\leq\omega$. In the final walk with $n=\omega$,
either the customer rents a bike at a nonempty station, or he directly leaves
the bike sharing system if no bike is rented after $\omega$\ consequent walks.
Thus the number of the consequent walks to find an available station may be
$1$ with probability $y_{0}^{\left(  N\right)  }\left(  t\right)  $, $2$ with
$\left[  y_{0}^{\left(  N\right)  }\left(  t\right)  \right]  ^{2}$, and
generally, $n$ with $\left[  y_{0}^{\left(  N\right)  }\left(  t\right)
\right]  ^{n}$ for $1\leq n\leq\omega$. Based on this, for the virtual
time-inhomogeneous $M(t)/M(t)/1/K$ queue, we obtain its service rates in
States $k$ for $1\leq k\leq K$ as follows:%
\begin{align*}
\eta_{k}^{\left(  N\right)  }\left(  t\right)   &  =\lambda+\gamma
y_{0}^{\left(  N\right)  }\left(  t\right)  +\gamma\left[  y_{0}^{\left(
N\right)  }\left(  t\right)  \right]  ^{2}+\gamma\left[  y_{0}^{\left(
N\right)  }\left(  t\right)  \right]  ^{3}+\cdots+\gamma\left[  y_{0}^{\left(
N\right)  }\left(  t\right)  \right]  ^{\omega}\\
&  =\lambda+\gamma y_{0}^{\left(  N\right)  }\left(  t\right)  \frac{1-\left[
y_{0}^{\left(  N\right)  }\left(  t\right)  \right]  ^{\omega}}{1-y_{0}%
^{\left(  N\right)  }\left(  t\right)  }\\
&  =\eta^{\left(  N\right)  }\left(  t\right)  ,
\end{align*}
which is independent of the number $k=1,2,\ldots,K$.

Now, we prove Equation (\ref{Equ-2}) in terms of the mean-field theory. Note
that we can compute the arrival rates $\xi_{l}^{\left(  N\right)  }\left(
t\right)  $ for $0\leq l\leq K-1$ according to a detailed probability analysis
on States $l$ for $0\leq l\leq K-1$.

For $l=0$ (i.e., States $0$), all the original $C$ bikes in the tagged station
are rented to travel on the roads. For the other $N-1$ stations, our
computation for the number of bikes rented to travel on the roads is based on
the mean-field theory (i.e., under an average setting), thus the expected
number of bikes rented to travel on the roads is given by%
\[
\left(  N-1\right)  \cdot\left[  C-\sum_{k=1}^{K}ky_{k}^{\left(  N\right)
}\left(  t\right)  \right]  ,
\]
where $\sum_{k=1}^{K}ky_{k}^{\left(  N\right)  }\left(  t\right)  $ is the
expected number of bikes parked in the tagged station, while $C-\sum_{k=1}%
^{K}ky_{k}^{\left(  N\right)  }\left(  t\right)  $ is the expected number of
bikes rented to travel on the roads from the tagged station. Therefore, for
the $N$ stations, the total expected number of bikes rented to travel on the
roads is given by
\[
C+\left(  N-1\right)  \cdot\left[  C-\sum_{k=1}^{K}ky_{k}^{\left(  N\right)
}\left(  t\right)  \right]  .
\]
Note that the returning-bike process of each bike is persistent in the sense
that the customer keeps finding an empty position in the next station, it is
easy to check that the return rate of each riding bike arriving at the tagged
station is given by%
\[
\mu+\mu y_{K}^{\left(  N\right)  }\left(  t\right)  +\mu\left[  y_{K}^{\left(
N\right)  }\left(  t\right)  \right]  ^{2}+\mu\left[  y_{K}^{\left(  N\right)
}\left(  t\right)  \right]  ^{3}+\cdots=\mu\frac{1}{1-y_{K}^{\left(  N\right)
}\left(  t\right)  },
\]
where $\left[  y_{K}^{\left(  N\right)  }\left(  t\right)  \right]  ^{n}$ is
the probability that a customer $n$ times continuously returns his bike to $n$
full stations. Thus we use the mean-field computation to obtain that for State
0 (for $l=0$),%
\begin{align*}
\xi_{0}^{\left(  N\right)  }\left(  t\right)   &  =\frac{1}{N}\left\{
C+\left(  N-1\right)  \left[  C-\sum_{k=1}^{K}ky_{k}^{\left(  N\right)
}\left(  t\right)  \right]  \right\}  \cdot\mu\frac{1}{1-y_{K}^{\left(
N\right)  }\left(  t\right)  }\\
&  =\frac{\mu}{N}\frac{1}{1-y_{K}^{\left(  N\right)  }\left(  t\right)
}\left\{  C+\left(  N-1\right)  \left[  C-\sum_{k=1}^{K}ky_{k}^{\left(
N\right)  }\left(  t\right)  \right]  \right\}  .
\end{align*}
Similarly, for States $l$ with $1\leq l\leq C-1$, we have%
\[
\xi_{l}^{\left(  N\right)  }\left(  t\right)  =\frac{\mu}{N}\frac{1}%
{1-y_{K}^{\left(  N\right)  }\left(  t\right)  }\left\{  \left(  C-l\right)
+\left(  N-1\right)  \left[  C-\sum_{k=1}^{K}ky_{k}^{\left(  N\right)
}\left(  t\right)  \right]  \right\}  .
\]
Finally, for States $l$ with $C\leq l\leq K$, since all the original $C$ bikes
are parked in the tagged station, we obtain%
\[
\xi_{l}^{\left(  N\right)  }\left(  t\right)  =\frac{\mu}{N}\frac{1}%
{1-y_{K}^{\left(  N\right)  }\left(  t\right)  }\left\{  \left(  N-1\right)
\left[  C-\sum_{k=1}^{K}ky_{k}^{\left(  N\right)  }\left(  t\right)  \right]
\right\}  ,
\]
which is independent of the number $l=C,C+1,\ldots,K$. This completes this
proof. \textbf{{\rule{0.08in}{0.08in}}}

\begin{Rem}
(1) In Equation (\ref{Equ-1}), for the number of consequent walks, it may be
useful to observe two special cases: (a) If $\omega=0$, $\eta_{k}^{\left(
N\right)  }\left(  t\right)  =\lambda$. (b) If $\omega\rightarrow\infty$, then
$\eta_{k}^{\left(  N\right)  }\left(  t\right)  =\lambda+\gamma y_{0}^{\left(
N\right)  }\left(  t\right)  /\left[  1-y_{0}^{\left(  N\right)  }\left(
t\right)  \right]  $.

(2) The time-inhomogeneous $M(t)/M(t)/1/K$ queue is a fictitious queueing
system corresponding to the number of bikes parked in the tagged station,
while its virtual arrival and virtual service rates are determined by means of
the empirical measure process through some mean-field computation.

(3) It is seen from the proof of Theorem 1 that the different ages of
``finding-bike attempts'' and ``returning-bike attempts'' has not any
influence on the mean-field computation due to the memoryless property of the
exponential distributions and of the Poisson processes. Thus, the mean-field
method can be successfully applied to our current analysis of the bike sharing
system. However, it will be very difficult (or an open problem) to apply the
mean-field method if there exist general distributions or general renewal
processes in the bike sharing system.
\end{Rem}

In the remainder of this section, we set up a system of mean-field equations
by means of the time-inhomogeneous $M(t)/M(t)/1/K$ queue whose state
transition relation is depcited in Figure 2 with the arrival rate $\xi
_{l}^{\left(  N\right)  }\left(  t\right)  $ for $0\leq l\leq K-1$, and with
service rate $\eta_{k}^{\left(  N\right)  }\left(  t\right)  =\eta^{\left(
N\right)  }\left(  t\right)  $ for $1\leq k\leq K$. To establish such
mean-field equations, readers nay refer to, such as, Li and Lui
\cite{Li:2014b}, Li et al. \cite{Li:2014a, Li:2015} and Fricker and Gast
\cite{Fri:2014} for more details.

To apply the mean-field theory, the number of bikes parked in the tagged
station is described as the virtual time-inhomogeneous $M(t)/M(t)/1/K$ queue.
Thus we can set up a system of mean-field equations in terms of the
(nonlinear) birth-death process corresponding to the $M(t)/M(t)/1/K$ queue. To
this end, we denote by $Q\left(  t\right)  $ the queue length of the
$M(t)/M(t)/1/K$ queue at time $t\geq0$. Then it is seen from Figure 2 that
$\left\{  Q\left(  t\right)  :t\geq0\right\}  $ is a time-inhomogeneous
continuous-time birth-death process whose infinite generator is given by%
\begin{equation}
\mathbf{V}_{\mathbf{y}^{\left(  N\right)  }\left(  t\right)  }=\left(
\begin{array}
[c]{cccccc}%
B_{1}\left(  t\right)  & B_{0}\left(  t\right)  &  &  &  & \\
B_{2}\left(  t\right)  & -\Theta_{C}^{\left(  N\right)  }\left(  t\right)  &
\xi_{C}^{\left(  N\right)  }\left(  t\right)  &  &  & \\
& \eta^{\left(  N\right)  }\left(  t\right)  & -\Theta_{C}^{\left(  N\right)
}\left(  t\right)  & \xi_{C}^{\left(  N\right)  }\left(  t\right)  &  & \\
&  & \ddots & \ddots & \ddots & \\
&  &  & \eta^{\left(  N\right)  }\left(  t\right)  & -\Theta_{C}^{\left(
N\right)  }\left(  t\right)  & \xi_{C}^{\left(  N\right)  }\left(  t\right) \\
&  &  &  & \eta^{\left(  N\right)  }\left(  t\right)  & -\eta^{\left(
N\right)  }\left(  t\right)
\end{array}
\right)  , \label{Equ0-0}%
\end{equation}
where%
\[
\eta^{\left(  N\right)  }\left(  t\right)  =\lambda+\gamma y_{0}^{\left(
N\right)  }\left(  t\right)  \frac{1-\left[  y_{0}^{\left(  N\right)  }\left(
t\right)  \right]  ^{\omega}}{1-y_{0}^{\left(  N\right)  }\left(  t\right)
},
\]
for $0\leq l\leq C$%
\[
\xi_{l}^{\left(  N\right)  }\left(  t\right)  =\frac{\mu}{N}\frac{1}%
{1-y_{K}^{\left(  N\right)  }\left(  t\right)  }\left\{  \left(  C-l\right)
+\left(  N-1\right)  \left[  C-\sum_{k=1}^{K}ky_{k}^{\left(  N\right)
}\left(  t\right)  \right]  \right\}
\]
and%
\[
\Theta_{l}^{\left(  N\right)  }\left(  t\right)  =\xi_{l}^{\left(  N\right)
}\left(  t\right)  +\eta^{\left(  N\right)  }\left(  t\right)  ;
\]%
\[
B_{1}\left(  t\right)  =\left(
\begin{array}
[c]{ccccc}%
-\xi_{0}^{\left(  N\right)  }\left(  t\right)  & \xi_{0}^{\left(  N\right)
}\left(  t\right)  &  &  & \\
\eta^{\left(  N\right)  }\left(  t\right)  & -\Theta_{1}^{\left(  N\right)
}\left(  t\right)  & \xi_{1}^{\left(  N\right)  }\left(  t\right)  &  & \\
& \ddots & \ddots & \ddots & \\
&  & \eta^{\left(  N\right)  }\left(  t\right)  & -\Theta_{C-2}^{\left(
N\right)  }\left(  t\right)  & \xi_{C-2}^{\left(  N\right)  }\left(  t\right)
\\
&  &  & \eta^{\left(  N\right)  }\left(  t\right)  & -\Theta_{C-1}^{\left(
N\right)  }\left(  t\right)
\end{array}
\right)  _{C\times C},
\]%
\[
B_{0}\left(  t\right)  =\left(  0,0,\ldots,0,\xi_{C-1}^{\left(  N\right)
}\left(  t\right)  \right)  ^{T}%
\]
and%
\[
B_{2}\left(  t\right)  =\left(  0,0,\ldots,0,\eta^{\left(  N\right)  }\left(
t\right)  \right)  ,
\]
$A^{T}$ denotes the transpose of the vector (or matrix) $A$.

Using the birth-death process described in Figure 2, we obtain a system of
mean-field (or ordinary differential) equations as follows:%
\[
\frac{\text{d}}{\text{d}t}y_{0}^{\left(  N\right)  }\left(  t\right)
=-\xi_{0}^{\left(  N\right)  }\left(  t\right)  y_{0}^{\left(  N\right)
}\left(  t\right)  +\eta^{\left(  N\right)  }\left(  t\right)  y_{1}^{\left(
N\right)  }\left(  t\right)  ,
\]
for $1\leq k\leq K-1$%
\[
\frac{\text{d}}{\text{d}t}y_{k}^{\left(  N\right)  }\left(  t\right)
=\xi_{k-1}^{\left(  N\right)  }\left(  t\right)  y_{k-1}^{\left(  N\right)
}\left(  t\right)  -\left[  \xi_{k}^{\left(  N\right)  }\left(  t\right)
+\eta^{\left(  N\right)  }\left(  t\right)  \right]  y_{k}^{\left(  N\right)
}\left(  t\right)  +\eta^{\left(  N\right)  }\left(  t\right)  y_{k+1}%
^{\left(  N\right)  }\left(  t\right)  ,
\]%
\[
\frac{\text{d}}{\text{d}t}y_{K}^{\left(  N\right)  }\left(  t\right)
=\xi_{K-1}^{\left(  N\right)  }\left(  t\right)  y_{K-1}^{\left(  N\right)
}\left(  t\right)  -\eta^{\left(  N\right)  }\left(  t\right)  y_{K}^{\left(
N\right)  }\left(  t\right)  .
\]

Now, we write the above system of mean-field equations in a vector form as%
\begin{equation}
\frac{\text{d}}{\text{d}t}\mathbf{y}^{\left(  N\right)  }\left(  t\right)
=\mathbf{y}^{\left(  N\right)  }\left(  t\right)  \mathbf{V}_{\mathbf{y}%
^{\left(  N\right)  }\left(  t\right)  }, \label{Equ0-1}%
\end{equation}
with the boundary and initial conditions%
\begin{equation}
\mathbf{y}^{\left(  N\right)  }\left(  t\right)  e=1,\text{ \ }\mathbf{y}%
^{\left(  N\right)  }\left(  0\right)  =\mathbf{g,} \label{Equ0-2}%
\end{equation}
where $\mathbf{g}=\left(  g_{0},g_{1},\ldots,g_{K}\right)  $ with $g_{i}\geq0$
for $0\leq i\leq K$ and $\sum_{i=0}^{K}g_{i}=1$, and $e$ is a column vector of
ones with a suitable dimension in the context.

\begin{Rem}
To deal with the time-inhomogeneous continuous-time birth-death process,
readers may refer to Chapter 8 in Li \cite{Li:2010} for more details, where
the detailed literatures are surveyed both for the time-inhomogeneous queues
and for the time-inhomogeneous Markov processes.
\end{Rem}

\section{A Lipschitz Condition}

In this section, we first establish a Lipschitz condition. Then we prove the
existence and uniqueness of solution to the system of ordinary differential
equations by means of the Lipschitz condition.

We write%
\begin{equation}
\frac{\text{d}}{\text{d}t}\mathbf{y}\left(  t\right)  =\mathbf{y}\left(
t\right)  \mathbf{V}_{\mathbf{y}\left(  t\right)  }, \label{Equa-1}%
\end{equation}
with the boundary and initial conditions%
\begin{equation}
\mathbf{y}\left(  t\right)  e=1,\text{ \ }\mathbf{y}\left(  0\right)
=\mathbf{g}, \label{Equa-2}%
\end{equation}
where
\begin{equation}
\mathbf{V}_{\mathbf{y}\left(  t\right)  }=\left(
\begin{array}
[c]{ccccc}%
-a\left(  t\right)  & a\left(  t\right)  &  &  & \\
b\left(  t\right)  & -c\left(  t\right)  & a\left(  t\right)  &  & \\
& \ddots & \ddots & \ddots & \\
&  & b\left(  t\right)  & -c\left(  t\right)  & a\left(  t\right) \\
&  &  & b\left(  t\right)  & -b\left(  t\right)
\end{array}
\right)  , \label{Equa-3}%
\end{equation}%
\[
b\left(  t\right)  =\lambda+\gamma y_{0}\left(  t\right)  \frac{1-\left[
y_{0}\left(  t\right)  \right]  ^{\omega}}{1-y_{0}\left(  t\right)  },
\]%
\[
a\left(  t\right)  =\mu\frac{1}{1-y_{K}\left(  t\right)  }\left[  C-\sum
_{k=1}^{K}ky_{k}\left(  t\right)  \right]
\]
and%
\[
c\left(  t\right)  =a\left(  t\right)  +b\left(  t\right)  .
\]
Obviously, that Equations (\ref{Equa-1}) and (\ref{Equa-2}) are a system of
first-order ordinary differential equations.

\begin{Rem}
Note that the system of ordinary differential equations (\ref{Equa-1}) and
(\ref{Equa-2}) is the limiting version of Equations (\ref{Equ0-1}) and
(\ref{Equ0-2}) as $N\rightarrow\infty$, while the existence of the limit
$\lim_{N\rightarrow\infty}\mathbf{y}^{\left(  N\right)  }\left(  t\right)
=\mathbf{y}\left(  t\right)  $ will be proved in the next section according to
the martingale limits and the weak convergence in the Skorohod space.
\end{Rem}

To discuss the existence and uniqueness of solution to the system of ordinary
differential equations (\ref{Equa-1}) and (\ref{Equa-2}), in what follows we
need to establish a Lipschitz condition by means of a computational method
given in Section 4.1 of Li et al. \cite{Li:2014a}.

For simplicity of description, we first suppress time $t$ from the vector
$\mathbf{y}\left(  t\right)  $ and its entries $y_{k}\left(  t\right)  $ for
$0\leq k\leq K$. Then we rewrite Equations (\ref{Equa-1}) and (\ref{Equa-2})
in a simple form as%
\begin{equation}
\frac{\text{d}}{\text{d}t}\mathbf{y}=F\left(  \mathbf{y}\right)  ,\text{
\ \ \ \ }\mathbf{y}e\mathbf{=1,y}\left(  0\right)  =\mathbf{g}, \label{Equa-4}%
\end{equation}
where%
\[
F\left(  \mathbf{y}\right)  =\mathbf{yV}_{\mathbf{y}}=\left(  y_{0}%
,y_{1},\ldots,y_{K}\right)  \left(
\begin{array}
[c]{ccccc}%
-a & a &  &  & \\
b & -c & a &  & \\
& \ddots & \ddots & \ddots & \\
&  & b & -c & a\\
&  &  & b & -b
\end{array}
\right)  ,
\]%
\[
b=\lambda+\frac{\gamma y_{0}\left(  1-y_{0}^{\omega}\right)  }{1-y_{0}},\text{
\ \ }a=\frac{\mu}{1-y_{K}}\left(  C-\sum_{k=1}^{K}ky_{k}\right)  ,
\]%
\[
c=\frac{\mu}{1-y_{K}}\left(  C-\sum_{k=1}^{K}ky_{k}\right)  +\left(
\lambda+\frac{\gamma y_{0}\left(  1-y_{0}^{\omega}\right)  }{1-y_{0}}\right)
.
\]
Let%
\[
F\left(  \mathbf{y}\right)  =\left(  F_{0}\left(  \mathbf{y}\right)
,F_{1}\left(  \mathbf{y}\right)  ,\ldots,F_{K-1}\left(  \mathbf{y}\right)
,F_{K}\left(  \mathbf{y}\right)  \right)  .
\]
Then for $k=0$%
\[
F_{0}\left(  \mathbf{y}\right)  =-y_{0}\frac{\mu}{1-y_{K}}\left(  C-\sum
_{k=1}^{K}ky_{k}\right)  +y_{1}\left[  \lambda+\frac{\gamma y_{0}\left(
1-y_{0}^{\omega}\right)  }{1-y_{0}}\right]  ,
\]
for $1\leq i\leq K-1$%
\[
F_{k}\left(  \mathbf{y}\right)  =\left(  y_{i-1}-y_{i}\right)  \frac{\mu
}{1-y_{K}}\left(  C-\sum_{k=1}^{K}ky_{k}\right)  +\left(  y_{i}-y_{i+1}%
\right)  \left[  \lambda+\frac{\gamma y_{0}\left(  1-y_{0}^{\omega}\right)
}{1-y_{0}}\right]
\]
and for $k=K$%
\[
F_{K}\left(  \mathbf{y}\right)  =y_{K-1}\frac{\mu}{1-y_{K}}\left(
C-\sum_{k=1}^{K}ky_{k}\right)  -y_{K}\left[  \lambda+\frac{\gamma y_{0}\left(
1-y_{0}^{\omega}\right)  }{1-y_{0}}\right]  .
\]

Now, we define the norms of a vector $\mathbf{x=}\left(  x_{0},x_{1}%
,\ldots,x_{K}\right)  $ and a matrix $A=\left(  a_{i,j}\right)  _{0\leq
i,j\leq K}$ as follows:%
\[
\left\|  \mathbf{x}\right\|  =\max_{0\leq i\leq K}\left\{  \left|
x_{i}\right|  \right\}
\]
and%
\[
\left\|  A\right\|  =\max_{0\leq j\leq K}\left\{  \sum_{i=0}^{K}\left|
a_{i,j}\right|  \right\}  .
\]
It is easy to check that%
\[
\left\|  \mathbf{x}A\right\|  \leq\left\|  \mathbf{x}\right\|  \left\|
A\right\|  .
\]

From (41) of Li et al. \cite{Li:2014a}, the matrix of partial derivatives of
the vector function $F\left(  \mathbf{y}\right)  $ of dimension $K+1$ is given
by%
\begin{equation}
DF\left(  \mathbf{y}\right)  =\left(
\begin{array}
[c]{cccc}%
\frac{\partial F_{0}\left(  \mathbf{y}\right)  }{\partial y_{0}} &
\frac{\partial F_{1}\left(  \mathbf{y}\right)  }{\partial y_{0}} & \cdots &
\frac{\partial F_{K}\left(  \mathbf{y}\right)  }{\partial y_{0}}\\
\frac{\partial F_{0}\left(  \mathbf{y}\right)  }{\partial y_{1}} &
\frac{\partial F_{1}\left(  \mathbf{y}\right)  }{\partial y_{1}} & \cdots &
\frac{\partial F_{K}\left(  \mathbf{y}\right)  }{\partial y_{1}}\\
\vdots & \vdots &  & \vdots\\
\frac{\partial F_{0}\left(  \mathbf{y}\right)  }{\partial y_{K}} &
\frac{\partial F_{1}\left(  \mathbf{y}\right)  }{\partial y_{K}} &  &
\frac{\partial F_{K}\left(  \mathbf{y}\right)  }{\partial y_{K}}%
\end{array}
\right)  . \label{Equa-5}%
\end{equation}

To establish the Lipschitz condition of the vector function $F\left(
\mathbf{y}\right)  $ of dimension $K+1$, it is seen from Lemma 5 of Li et al.
\cite{Li:2014a} that we need to provide an upper bound of the norm $\left\|
DF\left(  \mathbf{y}\right)  \right\|  $. To this end, it is necessary to
first give an assumption with respect to the two key numbers $y_{0}$ and
$y_{K}$ as follows:

\textbf{Assumption of Problematic Stations:} Let $\delta$ be a sufficiently
small positive number. We assume that $0\leq y_{0},y_{K}\leq1-\delta$.

Now, we provide some interpretation for practical rationality of the
Assumption of Problematic Stations. Firstly, the probability $y_{0}\left(
t\right)  +y_{K}\left(  t\right)  $ of problematic stations is always smaller
by means of some management mechanism or control methods (for example,
repositioning by trucks, price incentives, and applications of information
technologies), thus it is natural and rational to take the condition: $0\leq
y_{0},y_{K}\leq1-\delta$ in practice. Secondly, Theorem 5 in Section 6 will
further demonstrate from the steady-state viewpoint that $\lim_{t\rightarrow
+\infty}y_{0}\left(  t\right)  =p_{0}\leq1/2$ and $\lim_{t\rightarrow+\infty
}y_{K}\left(  t\right)  =p_{K}\leq1-\delta$. Finally, if $y_{0}\left(
t\right)  =1$, then $y_{k}\left(  t\right)  =0$ for $1\leq k\leq K$; while if
$y_{K}\left(  t\right)  =1$, then $y_{k}\left(  t\right)  =0$ for $0\leq k\leq
K-1$. Therefore, such a case with either $y_{0}\left(  t\right)  =1$ or
$y_{K}\left(  t\right)  =1$ will directly lead to the unavailability of the
bike sharing system.

\begin{The}
\label{The:Lip}(1) Under the Assumption of Problematic Stations, $\left\|
DF\left(  \mathbf{y}\right)  \right\|  \leq\mathbf{M}$, where%
\[
\mathbf{M}=2\lambda+\gamma\frac{\omega\left(  \omega+5\right)  }{2}+\frac{\mu
}{\delta}\left[  \left(  1+\frac{1}{\delta}\right)  C+\frac{K\left(
K+1\right)  }{2}\right]  \mathbf{.}%
\]

(2) The the vector function $F\left(  \mathbf{y}\right)  $ of dimension $K+1$
is continuous and also satisfies the Lipschitz condition for $\left(
t,\mathbf{y}\right)  \in\lbrack0,+\infty)\times\left\{  \lbrack0,1-\delta
]\times\left[  0,1\right]  ^{K-1}\times\lbrack0,1-\delta]\right\}  $.

(3) There exists a unique solution to the system of ordinary differential
equations $\frac{\text{d}}{\text{d}t}\mathbf{y}=F\left(  \mathbf{y}\right)
$,\ $\mathbf{y}e=1$ and $\mathbf{y}\left(  0\right)  =\mathbf{g}$ for $\left(
t,\mathbf{y}\right)  \in\lbrack0,+\infty)\times\left\{  \lbrack0,1-\delta
]\times\left[  0,1\right]  ^{K-1}\times\lbrack0,1-\delta]\right\}  $.
\end{The}

\textbf{Proof: \ }(1)\ It follows from (\ref{Equa-5}) that%
\[
\left\|  DF\left(  \mathbf{y}\right)  \right\|  =\max_{0\leq j\leq K}\left\{
\sum_{i=0}^{K}\left|  \frac{\partial F_{j}\left(  \mathbf{y}\right)
}{\partial y_{i}}\right|  \right\}  .
\]
It is easy to check that%
\[
\frac{\partial F_{0}\left(  \mathbf{y}\right)  }{\partial y_{0}}=-\frac{\mu
}{1-y_{K}}\left(  C-\sum_{k=1}^{K}ky_{k}\right)  +\gamma y_{1}\sum
_{k=1}^{\omega-1}ky_{0}^{k},
\]%
\[
\frac{\partial F_{0}\left(  \mathbf{y}\right)  }{\partial y_{1}}%
=y_{0}\frac{\mu}{1-y_{K}}+\lambda+\gamma y_{0}\sum_{k=0}^{\omega-1}y_{0}^{k},
\]
and for $2\leq i\leq K$%
\[
\frac{\partial F_{0}\left(  \mathbf{y}\right)  }{\partial y_{i}}%
=y_{0}\frac{i\mu}{1-y_{K}}.
\]
By using%
\[
\left|  y_{k}\right|  \leq1,0\leq k\leq K;\text{ \ }\frac{1}{1-y_{K}}%
\leq\frac{1}{\delta};\text{ \ }0\leq C-\sum_{k=1}^{K}ky_{k}\leq C,
\]
we obtain%
\[
\sum_{i=0}^{K}\left|  \frac{\partial F_{0}\left(  \mathbf{y}\right)
}{\partial y_{i}}\right|  \leq\lambda+\frac{\mu}{\delta}\left(
C+\frac{K\left(  K+1\right)  }{2}\right)  +\gamma\frac{\omega\left(
\omega+3\right)  }{2}.
\]
Similarly, we obtain that for $1\leq j\leq K-1$%
\[
\sum_{i=0}^{K}\left|  \frac{\partial F_{j}\left(  \mathbf{y}\right)
}{\partial y_{i}}\right|  \leq2\lambda+\frac{\mu}{\delta}\left(
2C+\frac{K\left(  K+1\right)  }{2}\right)  +\gamma\frac{\omega\left(
\omega+5\right)  }{2}%
\]
and%
\[
\sum_{i=0}^{K}\left|  \frac{\partial F_{K}\left(  \mathbf{y}\right)
}{\partial y_{i}}\right|  \leq\lambda+\frac{\mu}{\delta}\left[  \left(
1+\frac{1}{\delta}\right)  C+\frac{K\left(  K+1\right)  }{2}\right]
+\gamma\frac{\omega\left(  \omega+3\right)  }{2}.
\]
Let%
\[
\mathbf{M}=2\lambda+\gamma\frac{\omega\left(  \omega+5\right)  }{2}+\frac{\mu
}{\delta}\left[  \left(  1+\frac{1}{\delta}\right)  C+\frac{K\left(
K+1\right)  }{2}\right]  .
\]
Then%
\[
\left\|  DF\left(  \mathbf{y}\right)  \right\|  =\max_{0\leq j\leq K}\left\{
\sum_{i=0}^{K}\left|  \frac{\partial F_{j}\left(  \mathbf{y}\right)
}{\partial y_{i}}\right|  \right\}  \leq\mathbf{M}.
\]

(2) By means of Lemma 5 in Li et al. \cite{Li:2014a}, we obtain that for any
two vectors $\mathbf{x,y\in}[0,1-\delta]\times\left[  0,1\right]  ^{K-1}%
\times\lbrack0,1-\delta]$,
\[
\left\|  F\left(  \mathbf{x}\right)  -F\left(  \mathbf{y}\right)  \right\|
\leq\sup_{0\leq\widetilde{t}\leq1}\left\|  DF\left(  \mathbf{x+}\widetilde
{t}\left(  \mathbf{y-x}\right)  \right)  \right\|  \left\|  \mathbf{y-x}%
\right\|  \leq\mathbf{M}\left\|  \mathbf{y-x}\right\|  .
\]
This shows that $F\left(  \mathbf{y}\right)  $ is continuous and also
satisfies the Lipschitz condition for $\left(  t,\mathbf{y}\right)  \in
\lbrack0,+\infty)\times\left\{  \lbrack0,1-\delta]\times\left[  0,1\right]
^{K-1}\times\lbrack0,1-\delta]\right\}  $.

(3) Note that $F\left(  \mathbf{y}\right)  $ is continuous and also satisfies
the Lipschitz condition for $\left(  t,\mathbf{y}\right)  \in\lbrack
0,+\infty)\times\left\{  \lbrack0,1-\delta]\times\left[  0,1\right]
^{K-1}\times\lbrack0,1-\delta]\right\}  $, it follows from Chapter 1 of Hale
\cite{Hale:1980} that there exists a unique solution to the system of ordinary
differential equations $\frac{\text{d}}{\text{d}t}\mathbf{y}=F\left(
\mathbf{y}\right)  $,\ $\mathbf{y}e=1$ and $\mathbf{y}\left(  0\right)
=\mathbf{g}$ for $\left(  t,\mathbf{y}\right)  \in\lbrack0,+\infty
)\times\left\{  \lbrack0,1-\delta]\times\left[  0,1\right]  ^{K-1}%
\times\lbrack0,1-\delta]\right\}  $. This completes the proof.
\textbf{{\rule{0.08in}{0.08in}}}

In the remainder of this section, we set up a simple relation between the two
systems of ordinary differential equations (\ref{Equ0-1}) and (\ref{Equ0-2});
and (\ref{Equa-1}) and (\ref{Equa-2}) through a limiting assumption
$\lim_{N\rightarrow\infty}\mathbf{y}^{\left(  N\right)  }\left(  t\right)
=\mathbf{y}\left(  t\right)  $, the correctness of which will further be
proved in the next section. To this end, from Equation (\ref{Equ0-1}) we set%
\[
G\left(  \mathbf{y}^{\left(  N\right)  }\left(  t\right)  \right)
=\mathbf{y}^{\left(  N\right)  }\left(  t\right)  \mathbf{V}_{\mathbf{y}%
^{\left(  N\right)  }\left(  t\right)  }%
\]
or a simple form by suppressing $t$
\[
G\left(  \mathbf{y}^{\left(  N\right)  }\right)  =\mathbf{y}^{\left(
N\right)  }\mathbf{V}_{\mathbf{y}^{\left(  N\right)  }}.
\]
It follows from (\ref{Equ0-1}) and (\ref{Equ0-2}) that%
\[
\frac{\text{d}}{\text{d}t}\mathbf{y}^{\left(  N\right)  }=G\left(
\mathbf{y}^{\left(  N\right)  }\right)  ,\text{ \ }\mathbf{y}^{\left(
N\right)  }e=1.
\]

By using $\lim_{N\rightarrow\infty}\mathbf{y}^{\left(  N\right)  }\left(
t\right)  =\mathbf{y}\left(  t\right)  $, we obtain that for $0\leq k\leq K-1$%
\[
\lim_{N\rightarrow\infty}\xi_{k}^{\left(  N\right)  }\left(  t\right)
=a\left(  t\right)  ,
\]
and%
\[
\lim_{N\rightarrow\infty}\eta^{\left(  N\right)  }\left(  t\right)  =b\left(
t\right)  .
\]
Thus comparing the vector $G\left(  \mathbf{y}^{\left(  N\right)  }\right)  $
with the vector $F\left(  \mathbf{y}\right)  $, we obtain%
\[
\lim_{N\rightarrow\infty}G\left(  \mathbf{y}^{\left(  N\right)  }\right)
=F\left(  \mathbf{y}\right)  .
\]
Since%
\[
\frac{\text{d}}{\text{d}t}\left(  \lim_{N\rightarrow\infty}\mathbf{y}^{\left(
N\right)  }\left(  t\right)  \right)  =\frac{\text{d}}{\text{d}t}%
\mathbf{y=}F\left(  \mathbf{y}\right)
\]
and%
\[
\lim_{N\rightarrow\infty}\left(  \frac{\text{d}}{\text{d}t}\mathbf{y}^{\left(
N\right)  }\left(  t\right)  \right)  =\lim_{N\rightarrow\infty}G\left(
\mathbf{y}^{\left(  N\right)  }\right)  \mathbf{=}F\left(  \mathbf{y}\right)
,
\]
we obtain%
\[
\frac{\text{d}}{\text{d}t}\left(  \lim_{N\rightarrow\infty}\mathbf{y}^{\left(
N\right)  }\left(  t\right)  \right)  =\lim_{N\rightarrow\infty}\left(
\frac{\text{d}}{\text{d}t}\mathbf{y}^{\left(  N\right)  }\left(  t\right)
\right)  .
\]

\section{The Martingale Limit}

In this section, we provide a martingale limit (i.e., the weak convergence in
the Skorohod space) for the sequence of empirical measure Markov processes in
the bike sharing system.

We define a $\left(  K+1\right)  $-dimensional simplex%
\[
\mathcal{F}=\left\{  f=\left(  f_{0},f_{1},\ldots,f_{K-1},f_{K}\right)
:f_{k}\geq0\text{ and }\sum\limits_{k=0}^{K}f_{k}=1\right\}  ,
\]
and endow $\mathcal{F}$ with the metric%
\[
d\left(  x,y\right)  =\sup_{0\leq k\leq K}\frac{\left|  x_{k}-y_{k}\right|
}{k+1},\text{ \ }x,y\in\mathcal{F}.
\]
Obviously, $d\left(  x,y\right)  \leq1$ for $x,y\in\mathcal{F}$. Under the
metric, the space $\mathcal{F}$ is compact, complete and separable. Let
$D_{\mathcal{F}}[0,+\infty)$ be the space of right-continuous paths with left
limits in $\mathcal{F}$ endowed with the Skorohod metric. For the Skorohod
space and the weak convergence, readers may refer to Billingsley
\cite{Bil:1968} and Chapter 3 of Ethier and Kurtz \cite{Eth:1986} for more details.

For the the empirical measure $\mathbf{Y}^{\left(  N\right)  }\left(
t\right)  $, we write%
\[
\mathbf{W}\left(  \mathbf{Y}^{\left(  N\right)  }\left(  t\right)  \right)
=\left(
\begin{array}
[c]{cccccc}%
A_{1}^{\left(  N\right)  }\left(  t\right)  & A_{0}^{\left(  N\right)
}\left(  t\right)  &  &  &  & \\
A_{2}^{\left(  N\right)  }\left(  t\right)  & -\Gamma_{C}^{\left(  N\right)
}\left(  t\right)  & \alpha_{C}^{\left(  N\right)  }\left(  t\right)  &  &  &
\\
& \beta^{\left(  N\right)  }\left(  t\right)  & -\Gamma_{C}^{\left(  N\right)
}\left(  t\right)  & \alpha_{C}^{\left(  N\right)  }\left(  t\right)  &  & \\
&  & \ddots & \ddots & \ddots & \\
&  &  & \beta^{\left(  N\right)  }\left(  t\right)  & -\Gamma_{C}^{\left(
N\right)  }\left(  t\right)  & \alpha_{C}^{\left(  N\right)  }\left(  t\right)
\\
&  &  &  & \beta^{\left(  N\right)  }\left(  t\right)  & -\beta^{\left(
N\right)  }\left(  t\right)
\end{array}
\right)  ,
\]
where%
\[
\beta^{\left(  N\right)  }\left(  t\right)  =\lambda+\gamma Y_{0}^{\left(
N\right)  }\left(  t\right)  \frac{1-\left[  Y_{0}^{\left(  N\right)  }\left(
t\right)  \right]  ^{\omega}}{1-Y_{0}^{\left(  N\right)  }\left(  t\right)
},
\]
for $0\leq l\leq C$%
\[
\alpha_{l}^{\left(  N\right)  }\left(  t\right)  =\frac{\mu}{N}\frac{1}%
{1-Y_{K}^{\left(  N\right)  }\left(  t\right)  }\left\{  \left(  C-l\right)
+\left(  N-1\right)  \left[  C-\sum_{k=1}^{K}kY_{k}^{\left(  N\right)
}\left(  t\right)  \right]  \right\}
\]
and%
\[
\Gamma_{l}^{\left(  N\right)  }\left(  t\right)  =\alpha_{l}^{\left(
N\right)  }\left(  t\right)  +\beta^{\left(  N\right)  }\left(  t\right)  ;
\]%
\[
A_{1}^{\left(  N\right)  }\left(  t\right)  =\left(
\begin{array}
[c]{ccccc}%
-\alpha_{0}^{\left(  N\right)  }\left(  t\right)  & \alpha_{0}^{\left(
N\right)  }\left(  t\right)  &  &  & \\
\beta^{\left(  N\right)  }\left(  t\right)  & -\Gamma_{1}^{\left(  N\right)
}\left(  t\right)  & \alpha_{1}^{\left(  N\right)  }\left(  t\right)  &  & \\
& \ddots & \ddots & \ddots & \\
&  & \beta^{\left(  N\right)  }\left(  t\right)  & -\Gamma_{C-2}^{\left(
N\right)  }\left(  t\right)  & \alpha_{C-2}^{\left(  N\right)  }\left(
t\right) \\
&  &  & \beta^{\left(  N\right)  }\left(  t\right)  & -\Gamma_{C-1}^{\left(
N\right)  }\left(  t\right)
\end{array}
\right)  _{C\times C},
\]%
\[
A_{0}^{\left(  N\right)  }\left(  t\right)  =\left(  0,0,\ldots,0,\alpha
_{C-1}^{\left(  N\right)  }\left(  t\right)  \right)  ^{T}%
\]
and%
\[
A_{2}^{\left(  N\right)  }\left(  t\right)  =\left(  0,0,\ldots,0,\beta
^{\left(  N\right)  }\left(  t\right)  \right)  .
\]

For the sequence $\left\{  \mathbf{Y}^{(N)}(t),t\geq0\right\}  $ of the
empirical measure Markov processes, by means of a similar computation for
setting up the system of mean-field equations (\ref{Equ0-1}) and
(\ref{Equ0-2}), we can obtain a system of stochastic differential equations as
follows:%
\begin{equation}
\frac{\text{d}}{\text{d}t}\mathbf{Y}^{\left(  N\right)  }\left(  t\right)
=\mathbf{Y}^{\left(  N\right)  }\left(  t\right)  \mathbf{W}\left(
\mathbf{Y}^{\left(  N\right)  }\left(  t\right)  \right)  , \label{MartE-1}%
\end{equation}
with the boundary and initial conditions%
\begin{equation}
\mathbf{Y}^{\left(  N\right)  }\left(  t\right)  e=1,\text{ \ }\mathbf{Y}%
^{\left(  N\right)  }\left(  0\right)  =\mathbf{g}. \label{MartE-2}%
\end{equation}

For the random vector $\mathbf{Y}\left(  t\right)  =\left(  Y_{0}\left(
t\right)  ,Y_{1}\left(  t\right)  ,\ldots,Y_{K}\left(  t\right)  \right)  $,
we write
\begin{equation}
\mathbf{W}\left(  \mathbf{Y}\left(  t\right)  \right)  =\left(
\begin{array}
[c]{ccccc}%
-\alpha\left(  t\right)  & \alpha\left(  t\right)  &  &  & \\
\beta\left(  t\right)  & -\tau\left(  t\right)  & \alpha\left(  t\right)  &  &
\\
& \ddots & \ddots & \ddots & \\
&  & \beta\left(  t\right)  & -\tau\left(  t\right)  & \alpha\left(  t\right)
\\
&  &  & \beta\left(  t\right)  & -\beta\left(  t\right)
\end{array}
\right)  , \label{Generator}%
\end{equation}%
\[
\beta\left(  t\right)  =\lambda+\gamma Y_{0}\left(  t\right)  \frac{1-\left[
Y_{0}\left(  t\right)  \right]  ^{\omega}}{1-Y_{0}\left(  t\right)  },
\]%
\[
\alpha\left(  t\right)  =\mu\frac{1}{1-Y_{K}\left(  t\right)  }\left[
C-\sum_{k=1}^{K}kY_{k}\left(  t\right)  \right]
\]
and%
\[
\tau\left(  t\right)  =\alpha\left(  t\right)  +\beta\left(  t\right)  .
\]
Based on this, we write
\begin{equation}
\frac{\text{d}}{\text{d}t}\mathbf{Y}\left(  t\right)  =\mathbf{Y}\left(
t\right)  \mathbf{W}\left(  \mathbf{Y}\left(  t\right)  \right)  ,
\label{MartE-3-0}%
\end{equation}
with the boundary condition%
\begin{equation}
\mathbf{Y}\left(  t\right)  e=1 \label{MartE-3-1}%
\end{equation}
and the initial condition%
\begin{equation}
\mathbf{Y}\left(  0\right)  =\mathbf{g}. \label{MartE-3-2}%
\end{equation}
Using a similar analysis to that in Theorem \ref{The:Lip}, we can show that
there exists a unique solution to the system of stochastic differential
equations (\ref{MartE-3-0}) to (\ref{MartE-3-2}), where the Assumption of
Problematic Stations is also necessary.

The following lemma is useful for discussing the mean-field limit
$\mathbf{Y}\left(  t\right)  =\lim_{N\rightarrow\infty}\mathbf{Y}^{\left(
N\right)  }\left(  t\right)  $ for $t\geq0$.

\begin{Lem}
\label{Lem:Mart} For the sequence $\left\{  \mathbf{Y}^{(N)}(t),t\geq
0\right\}  $ of Markov processes,%
\begin{equation}
\mathbf{M}^{\left(  N\right)  }\left(  t\right)  =\mathbf{Y}^{\left(
N\right)  }\left(  t\right)  -\mathbf{Y}^{\left(  N\right)  }\left(  0\right)
-\int_{0}^{t}\left\{  \mathbf{Y}^{\left(  N\right)  }\left(  x\right)
\mathbf{W}\left(  \mathbf{Y}^{\left(  N\right)  }\left(  x\right)  \right)
\right\}  \text{d}x \label{MartE-5}%
\end{equation}
is a martingale with respect to $N\geq1$.
\end{Lem}

\textbf{Proof:} Note that the generator of the Markov process $\left\{
\mathbf{Y}^{(N)}(t),t\geq0\right\}  $ is given by the matrix $\mathbf{W}%
\left(  \mathbf{Y}^{\left(  N\right)  }\left(  t\right)  \right)  $, thus
using Dynkin's formula, e.g., see Equation (III.10.13) in Rogers and Williams
\cite{Rog:1994} or Page 162 in Ethier and Kurtz \cite{Eth:1986}, and it is
easy to check that $\mathbf{M}^{\left(  N\right)  }\left(  t\right)  $ is a
martingale with respect to $N\geq1$. This completes the proof.
\textbf{{\rule{0.08in}{0.08in}}}

The following theorem gives the mean-field limit of the sequence $\left\{
\mathbf{Y}^{(N)}(t),t\geq0\right\}  $ of Markov processes. Note that this
mean-field limit is a key to proving the asymptotic independence of the bike
sharing system.

\begin{The}
\label{The:Lim}If $\mathbf{Y}^{\left(  N\right)  }\left(  0\right)
$\ converges weakly to$\ \mathbf{Y}\left(  0\right)  \in\mathcal{F}$ as
$N\rightarrow\infty$, then $\left\{  \mathbf{Y}^{\left(  N\right)  }\left(
t\right)  ,N\geq1\right\}  $ converges weakly in $D_{\mathcal{F}}[0,+\infty)$
endowed with the Skorohod topology to the solution $\mathbf{Y}\left(
t\right)  $ to the system of stochastic differential equations
(\ref{MartE-3-0}) to (\ref{MartE-3-2}).
\end{The}

\textbf{Proof: }The proof can be completed by the following three steps.

\textbf{Step One: }The relative compactness of $\mathbf{Y}^{\left(  N\right)
}\left(  t\right)  $ in $D_{\mathcal{F}}[0,+\infty)$

Note that the space $\mathcal{F}$ is of dimension $K+1$, we use Paragraphs 8.6
to 8.9 of Chapter 3 of Ethier and Kurtz \cite{Eth:1986} (see Pages 137 to 139)
to prove the relative compactness of $\mathbf{Y}^{\left(  N\right)  }\left(
t\right)  $ in $D_{\mathcal{F}}[0,+\infty)$. To that end, we only need to
indicate three conditions given in Chapter 3 of Ethier and Kurtz
\cite{Eth:1986} as follows:

\begin{description}
\item[(a) EK7.7] For every $\varepsilon>0$ and rational $r\geq0$, there exists
a compact set $\Gamma_{\varepsilon,r}\in\mathcal{F}$ such that%
\[
\lim_{N\rightarrow\infty}\inf_{y\in\Gamma_{\varepsilon,r}}P\left\{  d\left(
\mathbf{Y}^{\left(  N\right)  }\left(  t\right)  ,y\right)  <\varepsilon
\right\}  \geq1-\varepsilon.
\]

\item[(b) EK8.37] For all $T>0$, there exists $\chi>0$, $D>0$ and $\tau>1$
such that for all $N\geq1$ and all $0\leq h\leq t\leq T+1$%
\[
E\left[  d^{\frac{\chi}{2}}\left(  \mathbf{Y}^{\left(  N\right)  }\left(
t+h\right)  ,\mathbf{Y}^{\left(  N\right)  }\left(  t\right)  \right)
d^{\frac{\chi}{2}}\left(  \mathbf{Y}^{\left(  N\right)  }\left(  t\right)
,\mathbf{Y}^{\left(  N\right)  }\left(  t-h\right)  \right)  \right]  \leq
Dh^{\tau}.
\]

\item[(c) EK8.30] For the above value $\chi>0$%
\[
\lim_{\delta\rightarrow0}\lim_{N\rightarrow\infty}\sup E\left[  d^{\chi
}\left(  \mathbf{Y}^{\left(  N\right)  }\left(  \delta\right)  ,\mathbf{Y}%
^{\left(  N\right)  }\left(  0\right)  \right)  \right]  =0.
\]
\end{description}

In what follows we prove each of the three conditions.

Firstly, we prove (a) EK7.7. Taking $\Gamma_{\varepsilon,r}=\mathcal{F}$, and
note that the space $\mathcal{F}$ is compact, this directly gives the proof of
(a) EK7.7 through a similar analysis to that in Theorem 7.2 of Chapter 3 of
Ethier and Kurtz \cite{Eth:1986} (see Pages 128 to 129).

Secondly, we prove (b) EK8.37. Let $\chi=2$. Then by using Remark 8.9 of
Chapter 3 of Ethier and Kurtz \cite{Eth:1986} (see Page 139), we obtain%
\begin{align}
&  E\left[  d^{\frac{\chi}{2}}\left(  \mathbf{Y}^{\left(  N\right)  }\left(
t+h\right)  ,\mathbf{Y}^{\left(  N\right)  }\left(  t\right)  \right)
d^{\frac{\chi}{2}}\left(  \mathbf{Y}^{\left(  N\right)  }\left(  t\right)
,\mathbf{Y}^{\left(  N\right)  }\left(  t-h\right)  \right)  \right]
\nonumber\\
&  =E\left[  d\left(  \mathbf{Y}^{\left(  N\right)  }\left(  t+h\right)
,\mathbf{Y}^{\left(  N\right)  }\left(  t\right)  \right)  \right]  \cdot
E\left[  d\left(  \mathbf{Y}^{\left(  N\right)  }\left(  t\right)
,\mathbf{Y}^{\left(  N\right)  }\left(  t-h\right)  \right)  \right]
\nonumber\\
&  \leq\left[  \left(  \lambda+\mu+\gamma\right)  h\right]  ^{2},
\label{MartL}%
\end{align}
this indicates that (b) EK8.37 holds for the parameters: $T$, $t$, $h$,
$D=\left(  \lambda+\mu+\gamma\right)  ^{2}$ and $\tau=2$.

Finally, we prove (c) EK8.30. It follows from (\ref{MartL}) that%
\[
E\left[  d^{\chi}\left(  \mathbf{Y}^{\left(  N\right)  }\left(  \delta\right)
,\mathbf{Y}^{\left(  N\right)  }\left(  0\right)  \right)  \right]
\leq\left[  \left(  \lambda+\mu+\gamma\right)  \delta\right]  ^{\chi},
\]
this gives%
\[
\lim_{\delta\rightarrow0}\lim_{N\rightarrow\infty}\sup E\left[  d^{\chi
}\left(  \mathbf{Y}^{\left(  N\right)  }\left(  \delta\right)  ,\mathbf{Y}%
^{\left(  N\right)  }\left(  0\right)  \right)  \right]  =0.
\]

\textbf{Step Two: }The weakly convergent limit of $\left\{  \mathbf{Y}%
^{\left(  N\right)  }\left(  t\right)  \right\}  $ has almost surely
continuous sample paths for $t\geq0$

For $\mathbf{Y}\in D_{\mathcal{F}}[0,+\infty)$, we define%
\[
J\left(  \mathbf{Y},u\right)  =\sup_{0\leq t\leq u}\left\{  d\left(
\mathbf{Y}\left(  t\right)  ,\mathbf{Y}\left(  t^{-}\right)  \right)
\right\}
\]
and%
\[
J\left(  \mathbf{Y}\right)  =\int_{0}^{+\infty}e^{-u}J\left(  \mathbf{Y}%
,u\right)  \text{d}u.
\]
Using Theorem 10.2 (a) of Chapter 3 of Ethier and Kurtz \cite{Eth:1986} (see
Page 148), it is easy to check that for all $N\geq1$ and $u\geq0$, $J\left(
\mathbf{Y}^{\left(  N\right)  },u\right)  \leq1/N$ almost surely, which leads
to $J\left(  \mathbf{Y}^{\left(  N\right)  }\right)  \leq1/N$ almost surely.
Thus, as $N\rightarrow\infty$, if $\mathbf{Y}^{\left(  N\right)  }\left(
t\right)  \Rightarrow\mathbf{Y}\left(  t\right)  $, then $\mathbf{Y}\left(
t\right)  $ is almost surely continuous if and only if $J\left(
\mathbf{Y}^{\left(  N\right)  }\right)  \Rightarrow0$, where ``$\Rightarrow$''
denotes the weak convergence.

\textbf{Step Three: }The martingale limit

Given the continuity of any limit point, using the continuous mapping theorem
(e.g., see Whitt \cite{Whi:2002}), we prove that Equations (\ref{MartE-3-0})
and (\ref{MartE-3-1}) are satisfied by any limit point: $\mathbf{Y}\left(
t\right)  =\lim_{N\rightarrow\infty}\mathbf{Y}^{\left(  N\right)  }\left(
t\right)  $ for\ $t\geq0$ as follows:

Using the martingale central limit theorem (e.g., see Theorem 1.4 of Chapter 7
of Ethier and Kurtz \cite{Eth:1986} in Page 339), it follows from Lemma
\ref{Lem:Mart} that as $N\rightarrow\infty$, $\langle M_{k}^{\left(  N\right)
}\left(  t\right)  \rangle\overset{P}{\rightarrow}0$ for $t\geq0$, where
$\langle\cdot\rangle$ denotes the quadratic variation. Note that $\langle
M_{k}^{\left(  N\right)  }\left(  t\right)  \rangle$ only changes at time $t$
when $M_{k}^{\left(  N\right)  }\left(  t\right)  $ jumps, and it increases by
the square of the jump sizes, while the jump sizes are of order $1/N$ and the
jump rates are of order $N$. Using a similar analysis to that in Theorem 2 of
Section 4, we can prove that there exists a unique solution to the system of
stochastic differential equations (\ref{MartE-3-0}) and (\ref{MartE-3-1}) for
any initial value. Noting the relative compactness of $\mathbf{Y}^{\left(
N\right)  }\left(  t\right)  $ in $D_{\mathcal{F}}[0,+\infty)$ and using
Chapter 3 of Ethier and Kurtz \cite{Eth:1986}, we prove that the sequence
$\left\{  \mathbf{Y}^{\left(  N\right)  }\left(  t\right)  ,N\geq1\right\}  $
of Markov processes converges in the space $D_{\mathcal{F}}[0,+\infty)$ to the
Markov process $\left\{  \mathbf{Y}\left(  t\right)  ,N\geq1\right\}  $. This
completes the proof. \textbf{{\rule{0.08in}{0.08in}}}

Finally, it is necessary to provide some interpretation on Theorem
\ref{The:Lim}. If $\lim_{N\rightarrow\infty}\mathbf{Y}^{\left(  N\right)
}\left(  0\right)  =\mathbf{y}\left(  0\right)  =\mathbf{g}\in$ $\Omega$ in
probability, then Theorem \ref{The:Lim} shows that $\mathbf{Y}\left(
t\right)  =\lim_{N\rightarrow\infty}\mathbf{Y}^{\left(  N\right)  }\left(
t\right)  $ is concentrated on the trajectory $\Im_{\mathbf{g}}=\left\{
\mathbf{y}(t,\mathbf{g}):t\geq0\right\}  $, where $\mathbf{y}(t,\mathbf{g}%
)=E\left[  \mathbf{Y}\left(  t\right)  \text{ }|\text{ }\mathbf{Y}\left(
0\right)  =\mathbf{g}\right]  $, and $\mathbf{y}(0,\mathbf{g})=\mathbf{g}$.
This indicates the functional strong law of large numbers for the time
evolution of the fraction of each state of this bike sharing system, thus the
sequence $\left\{  \mathbf{Y}^{\left(  N\right)  }\left(  t\right)
,t\geq0\right\}  $ of Markov processes converges weakly to the expected
fraction vector $\mathbf{y}(t,\mathbf{g})$ as $N\rightarrow\infty$, that is,
for any $T>0$%
\begin{equation}
\lim_{N\rightarrow\infty}\sup_{0\leq s\leq T}\left\|  \mathbf{Y}^{\left(
N\right)  }\left(  s\right)  -\mathbf{y}(s,\mathbf{g})\right\|  =0\text{ \ in
probability}. \label{EquK-16}%
\end{equation}

\begin{Rem}
To study the weak convergence in the the Skorohod space for the sequence
$\left\{  \mathbf{Y}^{\left(  N\right)  }\left(  t\right)  ,t\geq0\right\}  $
of Markov processes, there are three frequently used methods: (1) Operator
semigroups, e.g., Vvedenskaya et al. \cite{Vve:1996}, Li and Lui
\cite{Li:2014b}, and Li et al. \cite{Li:2014a, Li:2015}; (2) martingale
limits, for example, Turner \cite{Tur:1998}, and Graham \cite{Gra:2000,
Gra:2004}; and (3) density-dependent jump Markov processes, for instance,
Chapter 11 of Ethier and Kurtz \cite{Eth:1986}, and Mitzenmacher
\cite{Mit:1996}. Here, this paper takes the method of martingale limits to
establish an outline of such a proof.
\end{Rem}

\begin{Rem}
Under the weak convergence in the the Skorohod space for the sequence
$\{\mathbf{Y}^{\left(  N\right)  }\left(  t\right)  $, $t\geq0\}$ of Markov
processes, Theorem \ref{The:Lim} demonstrates the correctness of the system of
mean-field equations (\ref{Equa-1}) and (\ref{Equa-2}), i.e., as
$N\rightarrow\infty$, Equations (\ref{Equa-1}) and (\ref{Equa-2}) are the
limits of Equations (\ref{Equ0-1}) and (\ref{Equ0-2}), respectively.
\end{Rem}

\section{The Fixed Point and Nonlinear Analysis}

In this section, we analyze the fixed point of the limiting system of
mean-field equations. We first prove that the fixed point is unique in terms
of the Birkhoff center. Then we simply analyze the asymptotic independence of
the bike sharing system, and also discuss the limiting interchangeability with
respect to $N\rightarrow\infty$ and $t\rightarrow+\infty$. Note that the
uniqueness of the fixed point is a key in numerical computation of the fixed
point in terms of a system of nonlinear equations.

Let the vector $\mathbf{p}$ be the fixed point of the limiting expected
fraction vector $\mathbf{y}(t)$. Then%
\[
\mathbf{p}=\lim_{t\rightarrow+\infty}\mathbf{y}\left(  t\right)  ,
\]
where $\mathbf{p}=\left(  p_{0},p_{1},\ldots,p_{K-1},p_{K}\right)  $ and%
\[
p_{k}=\lim_{t\rightarrow+\infty}y_{k}\left(  t\right)  ,\text{ \ \ }0\leq
k\leq K.
\]
This gives%
\[
\mathbf{p}=\lim_{t\rightarrow+\infty}\lim_{N\rightarrow\infty}\mathbf{y}%
^{\left(  N\right)  }(t).
\]

We write%
\[
b\left(  \mathbf{p}\right)  =\lim_{t\rightarrow+\infty}b\left(  t\right)
=\lambda+\gamma p_{0}\frac{1-p_{0}^{\omega}}{1-p_{0}},
\]%
\[
a\left(  \mathbf{p}\right)  =\lim_{t\rightarrow+\infty}a\left(  t\right)
=\mu\frac{1}{1-p_{K}}\left(  C-\sum_{k=1}^{K}kp_{k}\right)
\]
and%
\[
c\left(  \mathbf{p}\right)  =a\left(  \mathbf{p}\right)  +b\left(
\mathbf{p}\right)  .
\]
Thus it follows from (\ref{Equa-3}) that%
\begin{equation}
\mathbf{V}_{\mathbf{p}}=\lim_{t\rightarrow+\infty}\mathbf{V}_{\mathbf{y}%
\left(  t\right)  }=\left(
\begin{array}
[c]{ccccc}%
-a\left(  \mathbf{p}\right)  & a\left(  \mathbf{p}\right)  &  &  & \\
b\left(  \mathbf{p}\right)  & -c\left(  \mathbf{p}\right)  & a\left(
\mathbf{p}\right)  &  & \\
& \ddots & \ddots & \ddots & \\
&  & b\left(  \mathbf{p}\right)  & -c\left(  \mathbf{p}\right)  & a\left(
\mathbf{p}\right) \\
&  &  & b\left(  \mathbf{p}\right)  & -b\left(  \mathbf{p}\right)
\end{array}
\right)  , \label{Sta-1}%
\end{equation}
which is the infinitesimal generator of an irreducible, aperiodic and
positive-recurrent birth-death process due to the fact that $a\left(
\mathbf{p}\right)  >0,b\left(  \mathbf{p}\right)  >0,$ and the size of the
matrix $\mathbf{V}_{\mathbf{p}}$ is finite.

On the other hand, it is easy to see that the matrix $\mathbf{V}%
_{\mathbf{y}\left(  t\right)  }$ given in (\ref{Generator}) is also the
infinitesimal generator of a continuous-time birth-death process with state
space $\left\{  0,1,\ldots,K\right\}  $. Since $a\left(  t\right)  >0$,
$b\left(  t\right)  >0$ and $\mathbf{V}_{\mathbf{y}\left(  t\right)  }e=0$,
the birth-death process $\mathbf{V}_{\mathbf{y}\left(  t\right)  }$ is
irreducible, aperiodic and positive-recurrent. In this case, it is seen from
Vvedenskaya et al. \cite{Vve:1996} or Mitzenmacher \cite{Mit:1996} that
\[
\lim_{t\rightarrow+\infty}\frac{\text{d}}{\text{d}t}\mathbf{y}\left(
t\right)  =0
\]
or%
\[
\lim_{t\rightarrow+\infty}\mathbf{y}\left(  t\right)  \mathbf{V}%
_{\mathbf{y}\left(  t\right)  }=0.
\]
Thus it follows from (\ref{Equa-1}) and (\ref{Equa-2}) that%
\begin{equation}
\left\{
\begin{array}
[c]{c}%
\mathbf{pV}_{\mathbf{p}}=0,\\
\mathbf{p}e=1.
\end{array}
\right.  \label{Sta-2}%
\end{equation}

\subsection{Expressions for the fixed point}

Note that the matrix $\mathbf{V}_{\mathbf{p}}$ may be viewed as the
infinitesimal generator of an irreducible, aperiodic and positive-recurrent
birth-death process who corresponds to the $M/M/1/K$ queue with arrival rate
$a\left(  \mathbf{p}\right)  $ and service rate $b\left(  \mathbf{p}\right)
$. Let $\rho\left(  \mathbf{p}\right)  =a\left(  \mathbf{p}\right)  /b\left(
\mathbf{p}\right)  $. It is easy to check that (a) if $\rho\left(
\mathbf{p}\right)  =1$, then%
\begin{equation}
p_{k}=\frac{1}{K+1},\text{ }0\leq k\leq K; \label{Sta-3-1}%
\end{equation}
and (b) if $\rho\left(  \mathbf{p}\right)  \neq1$, then%
\begin{equation}
p_{k}=\rho^{k}\left(  \mathbf{p}\right)  \frac{1-\rho\left(  \mathbf{p}%
\right)  }{1-\rho^{K+1}\left(  \mathbf{p}\right)  },\text{ }0\leq k\leq K.
\label{Sta-3}%
\end{equation}
This demonstrates that if $\rho\left(  \mathbf{p}\right)  \neq1$, then the
probability vector $\mathbf{p}$ is the fixed point of the following nonlinear
vector equation%
\begin{equation}
\mathbf{p=}\left(  \frac{1-\rho\left(  \mathbf{p}\right)  }{1-\rho
^{K+1}\left(  \mathbf{p}\right)  },\rho\left(  \mathbf{p}\right)
\frac{1-\rho\left(  \mathbf{p}\right)  }{1-\rho^{K+1}\left(  \mathbf{p}%
\right)  },\ldots,\rho^{K}\left(  \mathbf{p}\right)  \frac{1-\rho\left(
\mathbf{p}\right)  }{1-\rho^{K+1}\left(  \mathbf{p}\right)  }\right)  .
\label{Sta-4}%
\end{equation}
Note that Li \cite{Li:2015a} gave some iterative algorithms for computing the
fixed point $\mathbf{p}$ by means by the system of nonlinear equations
(\ref{Sta-2}) or (\ref{Sta-4}).

In the following, we set up another nonlinear vector equation satisfied by the
fixed point $\mathbf{p}$. Different from Equation (\ref{Sta-4}), the new
nonlinear vector equation can be employed to study a more general
block-structure bike sharing system with either a Markovian arrival process
(MAP) or a phase-type (PH) service time, e.g., see Li \cite{Li:2010} and Li
and Lui \cite{Li:2014b} for more details.

To solve the system of equations (\ref{Sta-2}) from a more general setting,
let $r_{\min}\left(  \mathbf{p}\right)  $ and $g_{\min}\left(  \mathbf{p}%
\right)  $ be the minimal nonnegative solutions to the following two nonlinear
equations%
\[
a\left(  \mathbf{p}\right)  -\left[  a\left(  \mathbf{p}\right)  +b\left(
\mathbf{p}\right)  \right]  r\left(  \mathbf{p}\right)  +b\left(
\mathbf{p}\right)  r^{2}\left(  \mathbf{p}\right)  =0
\]
and%
\[
a\left(  \mathbf{p}\right)  g^{2}\left(  \mathbf{p}\right)  -\left[  a\left(
\mathbf{p}\right)  +b\left(  \mathbf{p}\right)  \right]  g\left(
\mathbf{p}\right)  +b\left(  \mathbf{p}\right)  =0,
\]
respectively. Then%
\[
r_{\min}\left(  \mathbf{p}\right)  =\frac{a\left(  \mathbf{p}\right)
+b\left(  \mathbf{p}\right)  -\left|  a\left(  \mathbf{p}\right)  -b\left(
\mathbf{p}\right)  \right|  }{2b\left(  \mathbf{p}\right)  }%
\]
and%
\[
g_{\min}\left(  \mathbf{p}\right)  =\frac{a\left(  \mathbf{p}\right)
+b\left(  \mathbf{p}\right)  -\left|  a\left(  \mathbf{p}\right)  -b\left(
\mathbf{p}\right)  \right|  }{2a\left(  \mathbf{p}\right)  }.
\]
Clearly, we have%
\[
r_{\min}\left(  \mathbf{p}\right)  b\left(  \mathbf{p}\right)  =g_{\min
}\left(  \mathbf{p}\right)  a\left(  \mathbf{p}\right)  =\frac{a\left(
\mathbf{p}\right)  +b\left(  \mathbf{p}\right)  -\left|  a\left(
\mathbf{p}\right)  -b\left(  \mathbf{p}\right)  \right|  }{2}.
\]
Let%
\begin{align*}
\Omega_{\mathbf{p}} =  &  \left\{  \left(  r_{\min}\left(  \mathbf{p}\right)
,\frac{1}{r_{\min}\left(  \mathbf{p}\right)  }\right)  :a\left(
\mathbf{p}\right)  >b\left(  \mathbf{p}\right)  \right\} \\
&  \bigcup\left\{  \left(  \frac{1}{g_{\min}\left(  \mathbf{p}\right)
},g_{\min}\left(  \mathbf{p}\right)  \right)  :a\left(  \mathbf{p}\right)
<b\left(  \mathbf{p}\right)  \right\} \\
&  \bigcup\left\{  \left(  1,1\right)  :a\left(  \mathbf{p}\right)  =b\left(
\mathbf{p}\right)  \right\}  .
\end{align*}
Then for a pair $\left(  r\left(  \mathbf{p}\right)  ,g\left(  \mathbf{p}%
\right)  \right)  \in\Omega_{\mathbf{p}}$, we have%
\[
r\left(  \mathbf{p}\right)  g\left(  \mathbf{p}\right)  =1.
\]

The following theorem illustrates that each element of the fixed point
$\mathbf{p}$ is a combinational sum of two geometric solutions if $a\left(
\mathbf{p}\right)  \neq b\left(  \mathbf{p}\right)  $.

\begin{The}
\label{The:Fix}If $a\left(  \mathbf{p}\right)  \neq b\left(  \mathbf{p}%
\right)  $ and $\left(  r\left(  \mathbf{p}\right)  ,g\left(  \mathbf{p}%
\right)  \right)  \in\Omega_{\mathbf{p}}$, then for $0\leq k\leq K$,%
\begin{equation}
p_{k}=c_{1}r^{k}\left(  \mathbf{p}\right)  +c_{2}g^{K-k}\left(  \mathbf{p}%
\right)  , \label{Exp0}%
\end{equation}
where the two constants $c_{1}$ and $c_{2}$ are determined by%
\begin{equation}
\left\{
\begin{array}
[c]{c}%
c_{1}=\frac{\frac{g^{K-1}\left(  \mathbf{p}\right)  \left[  b\left(
\mathbf{p}\right)  -g\left(  \mathbf{p}\right)  a\left(  \mathbf{p}\right)
\right]  }{a\left(  \mathbf{p}\right)  -r\left(  \mathbf{p}\right)  b\left(
\mathbf{p}\right)  }}{\frac{g^{K-1}\left(  \mathbf{p}\right)  \left[  b\left(
\mathbf{p}\right)  -g\left(  \mathbf{p}\right)  a\left(  \mathbf{p}\right)
\right]  }{a\left(  \mathbf{p}\right)  -r\left(  \mathbf{p}\right)  b\left(
\mathbf{p}\right)  }\frac{1-r^{K+1}\left(  \mathbf{p}\right)  }{1-r\left(
\mathbf{p}\right)  }-\frac{1-g^{K+1}\left(  \mathbf{p}\right)  }{1-g\left(
\mathbf{p}\right)  }},\\
c_{2}=\frac{1}{\frac{g^{K-1}\left(  \mathbf{p}\right)  \left[  b\left(
\mathbf{p}\right)  -g\left(  \mathbf{p}\right)  a\left(  \mathbf{p}\right)
\right]  }{a\left(  \mathbf{p}\right)  -r\left(  \mathbf{p}\right)  b\left(
\mathbf{p}\right)  }\frac{1-r^{K+1}\left(  \mathbf{p}\right)  }{1-r\left(
\mathbf{p}\right)  }-\frac{1-g^{K+1}\left(  \mathbf{p}\right)  }{1-g\left(
\mathbf{p}\right)  }}.
\end{array}
\right.  \label{Exp1}%
\end{equation}
\end{The}

\textbf{Proof: \ }If $a\left(  \mathbf{p}\right)  \neq b\left(  \mathbf{p}%
\right)  $, then the proof contains three steps. Firstly, it is easy to check
that for $1\leq k\leq K-1$, $p_{k}=c_{1}r^{k}\left(  \mathbf{p}\right)
+c_{2}g^{K-k}\left(  \mathbf{p}\right)  $ with $\left(  r\left(
\mathbf{p}\right)  ,g\left(  \mathbf{p}\right)  \right)  \in\Omega
_{\mathbf{p}}$ can satisfy the equation%
\[
p_{k-1}a\left(  \mathbf{p}\right)  -p_{k}\left[  a\left(  \mathbf{p}\right)
+b\left(  \mathbf{p}\right)  \right]  +p_{k+1}b\left(  \mathbf{p}\right)  =0.
\]
Secondly, for $k=0,K$ we obtain%
\begin{equation}
-\left[  c_{1}+c_{2}g^{K}\left(  \mathbf{p}\right)  \right]  a\left(
\mathbf{p}\right)  +\left[  c_{1}r\left(  \mathbf{p}\right)  +c_{2}%
g^{K-1}\left(  \mathbf{p}\right)  \right]  b\left(  \mathbf{p}\right)  =0
\label{Exp2}%
\end{equation}
and%
\begin{equation}
\left[  c_{1}r^{K-1}\left(  \mathbf{p}\right)  +c_{2}g\left(  \mathbf{p}%
\right)  \right]  a\left(  \mathbf{p}\right)  -\left[  c_{1}r^{K}\left(
\mathbf{p}\right)  +c_{2}\right]  b\left(  \mathbf{p}\right)  =0. \label{Exp3}%
\end{equation}
It follows from (\ref{Exp2}) and (\ref{Exp3}) that%
\begin{equation}
c_{1}=\frac{g^{K-1}\left(  \mathbf{p}\right)  b\left(  \mathbf{p}\right)
-g^{K}\left(  \mathbf{p}\right)  a\left(  \mathbf{p}\right)  }{a\left(
\mathbf{p}\right)  -r\left(  \mathbf{p}\right)  b\left(  \mathbf{p}\right)
}c_{2} \label{Exp4}%
\end{equation}
and%
\begin{equation}
c_{1}=\frac{b\left(  \mathbf{p}\right)  -g\left(  \mathbf{p}\right)  a\left(
\mathbf{p}\right)  }{r^{K-1}\left(  \mathbf{p}\right)  a\left(  \mathbf{p}%
\right)  -r^{K}\left(  \mathbf{p}\right)  b\left(  \mathbf{p}\right)  }c_{2}
\label{Exp5}%
\end{equation}
respectively. Note that $r\left(  \mathbf{p}\right)  g\left(  \mathbf{p}%
\right)  =1$ for $\left(  r\left(  \mathbf{p}\right)  ,g\left(  \mathbf{p}%
\right)  \right)  \in\Omega_{\mathbf{p}}$, we have%
\begin{align*}
\frac{b\left(  \mathbf{p}\right)  -g\left(  \mathbf{p}\right)  a\left(
\mathbf{p}\right)  }{r^{K-1}\left(  \mathbf{p}\right)  a\left(  \mathbf{p}%
\right)  -r^{K}\left(  \mathbf{p}\right)  b\left(  \mathbf{p}\right)  }  &
=\frac{\frac{1}{r^{K-1}\left(  \mathbf{p}\right)  }\left[  b\left(
\mathbf{p}\right)  -g\left(  \mathbf{p}\right)  a\left(  \mathbf{p}\right)
\right]  }{a\left(  \mathbf{p}\right)  -r\left(  \mathbf{p}\right)  b\left(
\mathbf{p}\right)  }\\
&  =\frac{g^{K-1}\left(  \mathbf{p}\right)  b\left(  \mathbf{p}\right)
-g^{K}\left(  \mathbf{p}\right)  a\left(  \mathbf{p}\right)  }{a\left(
\mathbf{p}\right)  -r\left(  \mathbf{p}\right)  b\left(  \mathbf{p}\right)  },
\end{align*}
this demonstrates that (\ref{Exp4}) is the same as (\ref{Exp5}). Finally,
using (\ref{Exp0}) and $\sum_{k=0}^{K}p_{k}=1$ we obtain
\[
c_{1}\frac{1-r^{K+1}\left(  \mathbf{p}\right)  }{1-r\left(  \mathbf{p}\right)
}+c_{2}\frac{1-g^{K+1}\left(  \mathbf{p}\right)  }{1-g\left(  \mathbf{p}%
\right)  }=1,
\]
which, together with (\ref{Exp4}), follows (\ref{Exp1}) in order to express
the constants $c_{1}$ and $c_{2}$. This completes the proof.
\textbf{{\rule{0.08in}{0.08in}}}

Using Theorem \ref{The:Fix}, the probability vector $\mathbf{p}$ is the fixed
point of the following nonlinear vector equation%
\begin{equation}
\mathbf{p=}\left(  c_{1}+c_{2}g^{K}\left(  \mathbf{p}\right)  ,c_{1}r\left(
\mathbf{p}\right)  +c_{2}g^{K-1}\left(  \mathbf{p}\right)  ,\ldots
,c_{1}r^{K-1}\left(  \mathbf{p}\right)  +c_{2}g\left(  \mathbf{p}\right)
,c_{1}r^{K}\left(  \mathbf{p}\right)  +c_{2}\right)  . \label{Sta-6}%
\end{equation}

We write%
\[
\mathbb{S}_{\mathbf{p}}=\left\{  \mathbf{p}:\mathbf{pV}_{\mathbf{p}%
}=0,\mathbf{p}e=1\right\}  .
\]
Then it is clear that%
\begin{align*}
\mathbb{S}_{\mathbf{p}}  &  =\left\{  \mathbf{p}:p_{k}=\rho^{k}\left(
\mathbf{p}\right)  \frac{1-\rho\left(  \mathbf{p}\right)  }{1-\rho
^{K+1}\left(  \mathbf{p}\right)  },\text{ }0\leq k\leq K\right\} \\
&  =\left\{  \mathbf{p}:p_{k}=c_{1}r^{k}\left(  \mathbf{p}\right)
+c_{2}g^{K-k}\left(  \mathbf{p}\right)  ,\text{ }0\leq k\leq K\right\}  .
\end{align*}
Since the equation $\mathbf{pV}_{\mathbf{p}}=0$ (or $p_{k}=\rho^{k}\left(
\mathbf{p}\right)  \left[  1-\rho\left(  \mathbf{p}\right)  \right]  /\left[
1-\rho^{K+1}\left(  \mathbf{p}\right)  \right]  $, or $p_{k}=c_{1}r^{k}\left(
\mathbf{p}\right)  +c_{2}g^{K-k}\left(  \mathbf{p}\right)  ,0\leq k\leq K$) is
nonlinear, it is possible for a more complicated bike sharing system that
there are multiple elements (solutions) in the set $\mathbb{S}_{\mathbf{p}}$.
In fact, an argument by analytic function indicates that the elements of the
set $\mathbb{S}_{\mathbf{p}}$ are isolated.

To describe the isolated element structure of the set $\mathbb{S}_{\mathbf{p}%
}$, we often need to use the Birkhoff center of the mean-field dynamic system,
which leads to check whether the fixed point is unique or not.

\subsection{The Birkhoff center and uniqueness}

For the Birkhoff center, our discussion includes the following two cases:

\textbf{Case one: }$N\rightarrow\infty$. In this case, we denote a solution to
the system of differential equations (\ref{Equa-1}) and (\ref{Equa-2}) by
$\Phi\left(  t\right)  $. Thus, the Birkhoff center of the solution
$\Phi\left(  t\right)  $ is defined as%
\begin{align*}
\mathbf{\Theta}=  &  \left\{  \overline{P}\in\mathcal{F}:\overline{P}%
=\lim_{k\rightarrow\infty}\Phi\left(  t_{k}\right)  \text{ for any scale
sequence}\right. \\
&  \left.  \left\{  t_{k}\right\}  \text{ with }t_{l}\geq0\text{ for }%
l\geq1\text{ and }\lim_{k\rightarrow\infty}t_{k}=+\infty\right\}  .
\end{align*}
Note that perhaps $\mathbf{\Theta}$ contains the limit cycles or the
stationary points (i.e., the local extremum points or the saddle points), it
is clear that $\mathbb{S}_{\mathbf{p}}\subset\mathbf{\Theta}$. Obviously, the
limiting empirical measure Markov process $\left\{  \mathbf{Y}\left(
t\right)  :t\geq0\right\}  $ spends most of its time in the Birkhoff center
$\mathbf{\Theta}$.

\textbf{Case two: }$t\rightarrow+\infty$. In this case, we write%
\[
\pi^{\left(  N\right)  }=\lim_{t\rightarrow+\infty}\mathbf{y}^{\left(
N\right)  }\left(  t\right)  ,
\]
since for each $N=1,2,3,\ldots$, the bike sharing system with $N$ identical
stations is stable.

Let%
\begin{align*}
\Xi=  &  \left\{  \overline{\pi}\in\mathcal{F}:\overline{\pi}=\lim
_{k\rightarrow\infty}\pi^{\left(  N_{k}\right)  }\text{ for any positive
integer sequence}\right. \\
&  \left.  \left\{  N_{k}\right\}  \text{ with }1\leq N_{1}\leq N_{2}\leq
N_{3}\leq\cdots\text{ and }\lim_{k\rightarrow\infty}N_{k}=\infty\right\}  .
\end{align*}
It is easy to see that%
\[
\mathbb{S}_{\mathbf{p}}\subset\Xi\subset\mathbf{\Theta.}%
\]
Therefore, the set $\mathbf{\Theta}-\mathbb{S}_{\mathbf{p}}$ contains the
limit cycles or the saddle points.

Note that%
\[
\left\{
\begin{array}
[c]{c}%
\mathbf{pV}_{\mathbf{p}}=0,\\
\mathbf{p}e=1,
\end{array}
\right.
\]
this gives that for $k=0$%

\begin{equation}
-\mu p_{0}\left(  1-p_{0}\right)  \left(  C-\sum_{k=1}^{K}kp_{k}\right)
+p_{1}\left[  \lambda\left(  1-p_{0}\right)  +\gamma p_{0}\left(
1-p_{0}^{\omega}\right)  \right]  \left(  1-p_{K}\right)  =0, \label{Uniq-1}%
\end{equation}
for $1\leq$ $k\leq K-1$%

\begin{equation}
-\mu(1-p_{0})\left(  C-\sum_{k=1}^{K}kp_{k}\right)  (p_{k-1}-p_{k})+\left[
\lambda\left(  1-p_{0}\right)  +\gamma p_{0}\left(  1-p_{0}^{\omega}\right)
\right]  \left(  1-p_{K}\right)  (p_{k}-p_{k+1})=0, \label{Uniq-2}%
\end{equation}
and for $k=K$%

\begin{equation}
-\mu p_{K-1}(1-p_{0})\left(  C-\sum_{k=1}^{K}kp_{k}\right)  +p_{K}\left[
\lambda\left(  1-p_{0}\right)  +\gamma p_{0}\left(  1-p_{0}^{\omega}\right)
\right]  \left(  1-p_{K}\right)  =0, \label{Uniq-3}%
\end{equation}
with the boundary condition%

\begin{equation}
p_{0}+p_{1}+p_{2}+\cdots+p_{K}=1. \label{Uniq-4}%
\end{equation}
Note that under the Assumption of Problematic Stations (i.e. $0<p_{0}%
,p_{K}<1-\delta$), the system of nonlinear equations (\ref{Sta-2}) is the same
as the system of nonlinear equations (\ref{Uniq-1}) to (\ref{Uniq-4}).

The following theorem gives an important result: The fixed point
$\mathbf{p}\in\mathbb{S}_{\mathbf{p}}$ is unique. Notice that the uniqueness
of the fixed point plays a key role in numerical computation for performance
measures of the bike sharing system. On the other hand, this proof uses the
system of nonlinear equations (\ref{Uniq-1}) to (\ref{Uniq-4}) by means of the
fact that the two special solutions $\left(  1,0,\ldots,0,0\right)  $ and
$\left(  0,0,\ldots,0,1\right)  $ are not in the set $\mathbb{S}_{\mathbf{p}}$.

\begin{The}
\label{The:Uniq}Let $\left|  \mathbb{S}_{\mathbf{p}}\right|  $ denote the
number of elements of the set $\mathbb{S}_{\mathbf{p}}$. Then $\left|
\mathbb{S}_{\mathbf{p}}\right|  =1$. This shows that the fixed point is unique.
\end{The}

\textbf{Proof: }This proof has two parts: (1) The existence of the fixed point
$\mathbf{p}$, which is easily dealt with by the fact that $\mathbf{p}$ is the
stationary probability vector of the ergodic birth-death process
$\mathbf{V}_{\mathbf{p}}$; and (2) the uniqueness of the fixed point
$\mathbf{p}$, which can be proved by means of the unique point of intersection
either between the quadratic function $f_{0}\left(  p_{0}\right)  $ and the
polynomial function $h_{0}\left(  p_{0}\right)  $, or between the quadratic
function $f_{n}\left(  p_{n}\right)  $ and the linear function $h_{n}\left(
p_{n}\right)  $ for $1\leq n\leq K-1$ as follows.

Based on the system of nonlinear equations (\ref{Uniq-1}) to (\ref{Uniq-4}),
the uniqueness of the fixed point $\mathbf{p}$ is proved through the following
three steps:

\textbf{Step one: }Analyzing $p_{0}$. In this case, we write%
\[
f_{0}\left(  p_{0}\right)  =\mu p_{0}\left(  1-p_{0}\right)  \left(
C-\sum_{k=1}^{K}kp_{k}\right)
\]
and%
\[
h_{0}\left(  p_{0}\right)  =p_{1}\left[  \lambda\left(  1-p_{0}\right)
+\gamma p_{0}\left(  1-p_{0}^{\omega}\right)  \right]  \left(  1-p_{K}\right)
.
\]
It is easy to check that%
\[
f_{0}\left(  0\right)  =0,\text{ }f_{0}\left(  1\right)  =0,\text{ }%
f_{0}\left(  \frac{1}{2}\right)  =\frac{1}{4}\mu\left(  C-\sum_{k=1}^{K}%
kp_{k}\right)  >0,
\]
and for $p_{0}\in\left(  0,1\right)  $%
\begin{align*}
\frac{\text{d}}{\text{d}p_{0}}f_{0}\left(  p_{0}\right)   &  =\left(
1-2p_{0}\right)  \mu\left(  C-\sum_{k=1}^{K}kp_{k}\right) \\
&  =\left\{
\begin{array}
[c]{cc}%
>0, & 0<p_{0}<\frac{1}{2},\\
=0, & p_{0}=\frac{1}{2},\\
<0, & \frac{1}{2}<p_{0}<1,
\end{array}
\right.
\end{align*}
and%
\[
\frac{\text{d}^{2}}{\text{d}\left(  p_{0}\right)  ^{2}}f_{0}\left(
p_{0}\right)  =-2\mu\left(  C-\sum_{k=1}^{K}kp_{k}\right)  <0,
\]
this demonstrates that $f_{0}\left(  p_{0}\right)  $ is a concave function
with the maximal value $f_{0}\left(  \frac{1}{2}\right)  >0$ at $p_{0}=1/2$.

Now, we analyze the polynomial function $h_{0}\left(  p_{0}\right)  $ for
$p_{0}\in\left(  0,1\right)  $. It is easy to see that%
\[
h_{0}\left(  0\right)  =\lambda p_{1}\left(  1-p_{K}\right)  >0,\text{
\ }h_{0}\left(  1\right)  =0\text{.}%
\]
For $p_{0}\in\left(  0,1\right)  $%
\begin{align}
\frac{\text{d}}{\text{d}p_{0}}h_{0}\left(  p_{0}\right)   &  =\left[
\gamma-\lambda-\gamma\left(  1+\omega\right)  p_{0}^{\omega}\right]
p_{1}\left(  1-p_{K}\right) \nonumber\\
&  =\left\{
\begin{array}
[c]{cc}%
>0, & p_{0}>\sqrt[\omega]{\frac{\gamma-\lambda}{\gamma\left(  1+\omega\right)
}},\\
=0, & p_{0}=\sqrt[\omega]{\frac{\gamma-\lambda}{\gamma\left(  1+\omega\right)
}},\\
<0, & p_{0}<\sqrt[\omega]{\frac{\gamma-\lambda}{\gamma\left(  1+\omega\right)
}}.
\end{array}
\right.  \label{FixD-1}%
\end{align}
Since $h_{0}\left(  0\right)  >0$ and $h_{0}\left(  1\right)  =0$, it is seen
from (\ref{FixD-1}) that only one case: $p_{0}<\sqrt[\omega]{\frac{\gamma
-\lambda}{\gamma\left(  1+\omega\right)  }}$ can hold; while the other two
cases are incorrect because the derivative $\frac{\text{d}}{\text{d}p_{0}%
}h_{0}\left(  p_{0}\right)  \geq0$ for $p_{0}\in\left(  0,1\right)  $ can not
result in such two values: $h_{0}\left(  0\right)  >0$ and $h_{0}\left(
1\right)  =0$. Thus we obtain%
\[
p_{0}<\sqrt[\omega]{\frac{\gamma-\lambda}{\gamma\left(  1+\omega\right)  }%
}<\sqrt[\omega]{\frac{1}{\left(  1+\omega\right)  }}\leq1.
\]
Note that for $p_{0}\in\left(  0,1\right)  $%
\[
\frac{\text{d}^{2}}{\text{d}\left(  p_{0}\right)  ^{2}}h_{0}\left(
p_{0}\right)  =-\gamma\omega\left(  1+\omega\right)  p_{0}^{\omega-1}%
p_{1}\left(  1-p_{K}\right)  <0,
\]
thus $h_{0}\left(  p_{0}\right)  $ is a decreasing and concave function from
Point $\left(  0,h_{0}\left(  0\right)  \right)  $ to $\left(  1,0\right)  $
without any extreme value.

Based on the above analysis, it is seen from Figure 4 (a) that there exists a
unique solution to the nonlinear equation $f_{0}\left(  p_{0}\right)  =$
$h_{0}\left(  p_{0}\right)  $ for $p_{0}\in\left(  0,1-\delta\right)  $.

\begin{figure}[ptb]
\centering     \includegraphics[width=12cm]{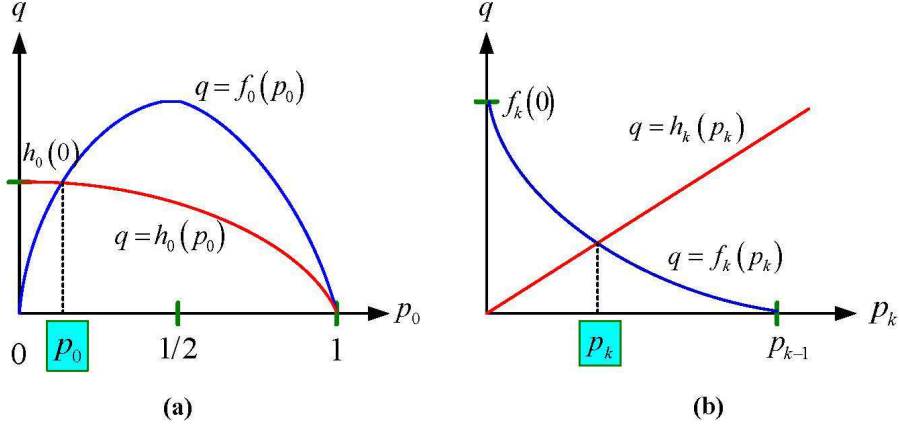}  \newline \caption{The
uniqueness of the fixed point}%
\label{figure:figure-4}%
\end{figure}

\textbf{Step two: }Analyzing $p_{k}$ for $1\leq k\leq K-1$. In this case, we
write%
\[
f_{k}\left(  p_{k}\right)  =\mu(1-p_{0})\left(  C-\sum_{k=1}^{K}kp_{k}\right)
(p_{k-1}-p_{k})
\]
and%
\[
h_{k}\left(  p_{k}\right)  =\left[  \lambda\left(  1-p_{0}\right)  +\gamma
p_{0}\left(  1-p_{0}^{\omega}\right)  \right]  \left(  1-p_{K}\right)
(p_{k}-p_{k+1}).
\]

Note that%
\[
f_{k}\left(  0\right)  =\mu(1-p_{0})\left(  C-\sum_{i\neq k}^{K}ip_{i}\right)
p_{k-1}>0,
\]%
\[
f_{k}\left(  p_{k-1}\right)  =0;
\]
and for $0<p_{k}<p_{k-1}$%
\[
\frac{\text{d}}{\text{d}p_{k}}f_{k}\left(  p_{k}\right)  =\mu(1-p_{0})\left[
-k(p_{k-1}-p_{k})-\left(  C-\sum_{k=1}^{K}kp_{k}\right)  \right]  <0,
\]%
\[
\frac{\text{d}^{2}}{\text{d}p_{k}^{2}}f_{k}\left(  p_{k}\right)
=2k\mu(1-p_{0})>0,
\]
thus the quadratic function $f_{k}\left(  p_{k}\right)  $ is a strictly
decreasing convex function for $0<p_{k}<p_{k-1}$.

Now, we consider the linear function $h_{k}\left(  p_{k}\right)  $. We obtain%
\[
h_{k}\left(  0\right)  =-\left[  \lambda\left(  1-p_{0}\right)  +\gamma
p_{0}\left(  1-p_{0}^{\omega}\right)  \right]  \left(  1-p_{K}\right)
p_{k+1}<0,
\]
and if $p_{k}=1$, then $p_{i}=0$ for $i\neq k$ with $1\leq i\leq K$, and it is
clear that%
\[
h_{k}\left(  1\right)  =\lambda>0.
\]
Since%
\[
\frac{\text{d}}{\text{d}p_{k}}h_{k}\left(  p_{k}\right)  =\left[
\lambda\left(  1-p_{0}\right)  +\gamma p_{0}\left(  1-p_{0}^{\omega}\right)
\right]  \left(  1-p_{K}\right)  >0,
\]
the linear function $h_{k}\left(  p_{k}\right)  $ is strictly increasing for
$p_{k}\in\left(  0,1\right)  $. Therefore, it is seen from Figure 4 (b) that
there exists a unique solution $p_{k}$ to the equation $f_{k}\left(
p_{k}\right)  =h_{k}\left(  p_{k}\right)  $.

\textbf{Step three: }Analyzing $p_{K}$. Since $p_{k}$ is the unique solution
to the equation $f_{k}\left(  p_{k}\right)  =h_{k}\left(  p_{k}\right)  $ for
$0\leq k\leq K-1$, it is clear that $p_{K}$ can uniquely determined by means
of the relation that $p_{K}=1-\sum_{k=0}^{K-1}p_{k}$. This completes the
proof. \textbf{{\rule{0.08in}{0.08in}}}

Now, we provide a simple discussion for the limiting interchangeability of the
vector $\mathbf{y}^{\left(  N\right)  }(t)$ as $N\rightarrow\infty$ and
$t\rightarrow+\infty$. Note that the limiting interchangeability is always
necessary and useful in many practical applications when using the stationary
probabilities of the limiting process $\left\{  \mathbf{Y}(t):t\geq0\right\}
$ to give an effective approximation for performance analysis of the bike
sharing system.

From $\left|  \mathbb{S}_{\mathbf{p}}\right|  =1$ by Theorem \ref{The:Uniq},
it is easy to see that%
\[
\lim_{t\rightarrow+\infty}\lim_{N\rightarrow\infty}\mathbf{y}^{\left(
N\right)  }\left(  t\right)  =\lim_{t\rightarrow+\infty}\mathbf{y}\left(
t\right)  =\mathbf{P}%
\]
and%
\[
\lim_{N\rightarrow\infty}\lim_{t\rightarrow+\infty}\mathbf{y}^{\left(
N\right)  }\left(  t\right)  =\lim_{N\rightarrow\infty}\mathbf{P}^{\left(
N\right)  }=\mathbf{P}.
\]
This gives%
\[
\lim_{t\rightarrow+\infty}\lim_{N\rightarrow\infty}\mathbf{y}^{\left(
N\right)  }(t)=\lim_{N\rightarrow\infty}\lim_{t\rightarrow+\infty}%
\mathbf{y}^{\left(  N\right)  }(t)=\mathbf{p}.
\]
Therefore, we have%
\[
\lim_{\substack{N\rightarrow\infty\\t\rightarrow+\infty}}\mathbf{y}^{\left(
N\right)  }(t)=\mathbf{p}.
\]

Finally, we provide a simple discussion on the asymptotic independence of this
bike sharing system. To this end, the uniqueness of the fixed point given by
$\left|  \mathbb{S}_{\mathbf{p}}\right|  =1$ of Theorem \ref{The:Uniq} plays a
key role. Using Corollaries 3 and 4 of Benaim and Le Boudec \cite{Ben:2008},
we obtain the asymptotic independence of the queueing processes of the bike
sharing system as follows:%
\begin{align*}
&  \lim_{t\rightarrow+\infty}\lim_{N\rightarrow\infty}P\left\{  X_{1}^{\left(
N\right)  }\left(  t\right)  =i_{1},X_{2}^{\left(  N\right)  }\left(
t\right)  =i_{2},\ldots,X_{k}^{\left(  N\right)  }\left(  t\right)
=i_{k}\right\} \\
&  =\lim_{N\rightarrow\infty}\lim_{t\rightarrow+\infty}P\left\{
X_{1}^{\left(  N\right)  }\left(  t\right)  =i_{1},X_{2}^{\left(  N\right)
}\left(  t\right)  =i_{2},\ldots,X_{k}^{\left(  N\right)  }\left(  t\right)
=i_{k}\right\} \\
&  =p_{i_{1}}p_{i_{2}}\cdots p_{i_{k}}%
\end{align*}
and%
\begin{align*}
&  \lim_{N\rightarrow\infty}\lim_{t\rightarrow+\infty}\frac{1}{t}\int_{0}%
^{t}\mathbf{1}_{\left\{  X_{1}^{\left(  N\right)  }\left(  t\right)
=i_{1},X_{2}^{\left(  N\right)  }\left(  t\right)  =i_{2},\ldots
,X_{k}^{\left(  N\right)  }\left(  t\right)  =i_{k}\right\}  }\text{d}t\\
&  =\lim_{t\rightarrow+\infty}\lim_{N\rightarrow\infty}\frac{1}{t}\int_{0}%
^{t}\mathbf{1}_{\left\{  X_{1}^{\left(  N\right)  }\left(  t\right)
=i_{1},X_{2}^{\left(  N\right)  }\left(  t\right)  =i_{2},\ldots
,X_{k}^{\left(  N\right)  }\left(  t\right)  =i_{k}\right\}  }\text{d}t\\
&  =p_{i_{1}}p_{i_{2}}\cdots p_{i_{k}}\text{ \ \ a.s.}%
\end{align*}

\begin{Rem}
For a more complicated bike sharing system, it is possible to have $\left|
\mathbb{S}_{\mathbf{p}}\right|  \geq2$. For this case with $\left|
\mathbb{S}_{\mathbf{p}}\right|  \geq2$, the metastability of the bike sharing
system is a key, and it can be roughly described as an interesting phenomenon
which occurs when the bike sharing system stays a very long time in some
abnormal state before reaching its normal state. To study the metastability, a
useful method is to determine a Lyapunov function $g\left(  \mathbf{y}\right)
$ for the system of differential equations (such as, (\ref{Equa-1}) and
(\ref{Equa-2})). Therefore, we need to find a continuously differentiable,
bounded from below, function $g\left(  \mathbf{y}\right)  $ defined on
$\left[  0,1\right]  ^{K+1}$ such that%
\[
\mathbf{yV}_{\mathbf{y}}\nabla g\left(  \mathbf{y}\right)  \leq0.
\]
Note that $\mathbf{yV}_{\mathbf{y}}\nabla g\left(  \mathbf{y}\right)  =0$ if
$\mathbf{yV}_{\mathbf{y}}=0$, which is satisfied by $\mathbf{y=p}$. On the
other hand, some properties of the function $g\left(  \mathbf{y}\right)  $
allow one to discriminate the stable points (the local minima of $g\left(
\mathbf{y}\right)  $) from the unstable points (the local maxima or saddle
points of $g\left(  \mathbf{y}\right)  $) in the study of metastability.
\end{Rem}

In general, it is not easy to give an analytic solution to the system of
nonlinear equations (\ref{Sta-2}), but its numerical solution may always be
simple and available. In the rest of this paper, we shall develop such a
numerical solution, and give numerical computation for performance measures of
this bike sharing system including the steady-state probability of the
problematic stations, and the stationary expected number of bikes at the
tagged station.

\section{Numerical Analysis}

In this section, we use some numerical examples to investigate the
steady-state probability of the problematic stations. Based on this,
performance analysis of the bike sharing system will focus on five points: (1)
$p_{0}$; (2) $p_{K}$; (3) $p_{0}+p_{K}$; (4) $E\left[  Q\right]  =\sum
_{k=1}^{K}kp_{k}$; and (5) the profit $R$.

Note that%
\[
\left\{
\begin{array}
[c]{c}%
\mathbf{pV}_{\mathbf{p}}=0,\\
\mathbf{p}e=1,
\end{array}
\right.
\]
this gives the system of nonlinear equations (\ref{Uniq-1}) to (\ref{Uniq-4})
whose solution is unique by means of $\left|  \mathbb{S}_{\mathbf{p}}\right|
=1$ by Theorem \ref{The:Uniq}. Also, we can numerically compute the unique
solution, i.e., the fixed point $\mathbf{p}$. Furthermore, the fixed point
$\mathbf{p}$ is employed in numerical computation for performance measures of
the bike sharing system. Based on this, we use some numerical examples to give
valuable observation and understanding with respect to design, operations and
optimization of the bike sharing systems. Therefore, such a numerical analysis
will become more and more useful in the study of bike sharing systems in practice.

\subsection{Analysis of $p_{0}$}

Note that $p_{0}$ is a probability that there is no bike in a tagged station,
thus it is also the probability that the arriving customer can not rent a bike
in the tagged station. To design a better bike sharing system, we hope that
the value of $p_{0}$ is as small as possible, and this can be realized through
taking a suitable parameters: $C,K,\lambda,\mu,\gamma$ and $\omega$, where
$C,K\ $and $\mu$ are controlled by the station; while $\lambda,\gamma$ and
$\omega$ are given by the customers.

In this bike sharing system, we take that $C=30$, $K=50$, $\omega=1$ and
$\gamma=0.25$. The left one of Figure 5 shows how the probability $p_{0}$
depends on $\lambda\in\left(  10,30\right)  $ when $\mu=0.3,1$ and $8$,
respectively. It is seen that $p_{0}$ increases either as $\lambda$ increases
or as $\mu$ decreases. Note that the numerical results are intuitively
reasonable because what $\lambda$ increases quickens up the rental rate of
bikes at the tagged station, while what $\mu$ decreases reduces the return
rate of bikes at the tagged station. Hence the probability $p_{0}$ increases
as the number of bikes parked at the tagged station decreases for the two cases.

\begin{figure}[ptb]
\centering              \includegraphics[width=7cm]{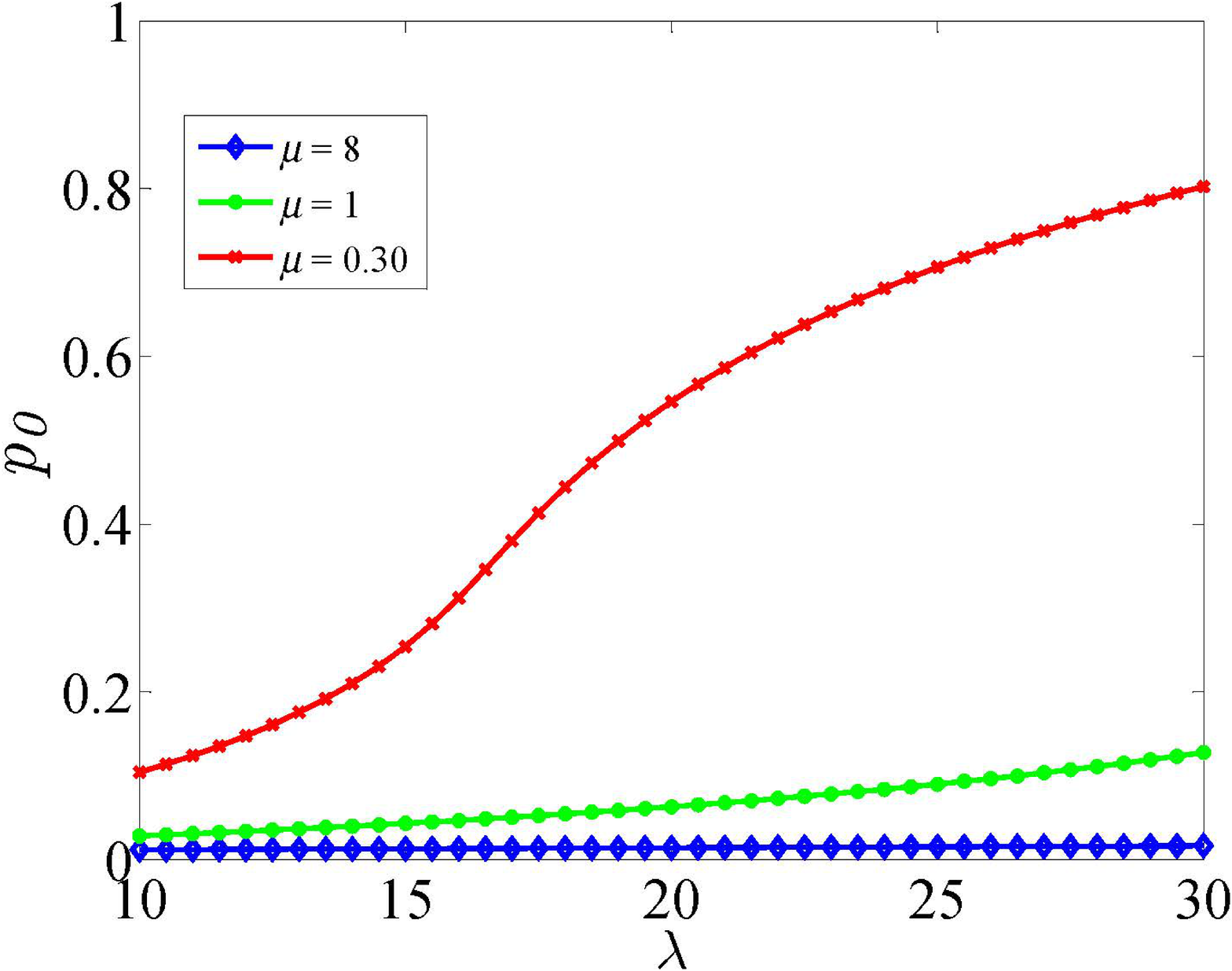}  \centering
\includegraphics[width=7cm]{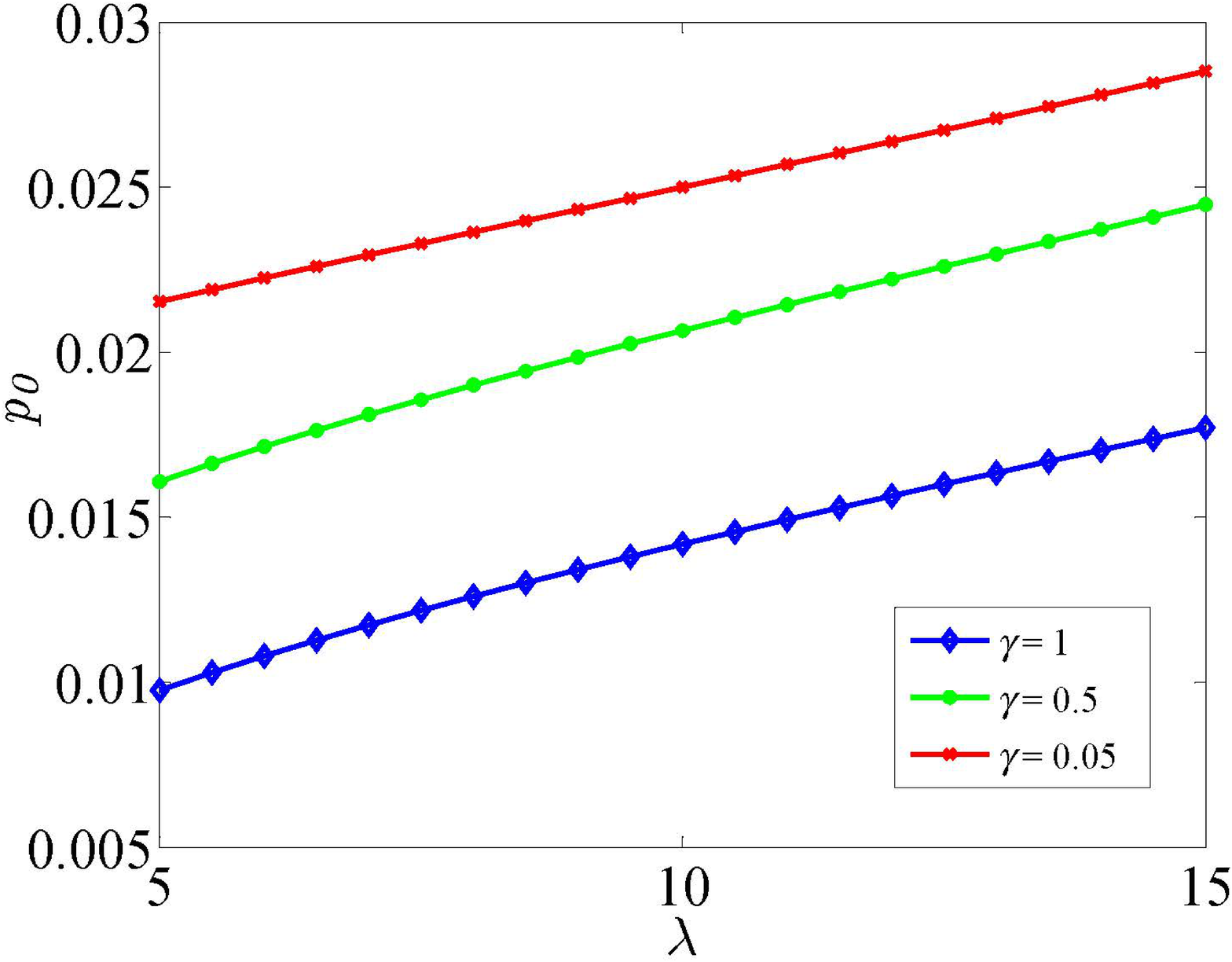}  \newline \caption{$p_{0}$ vs.
$\lambda$, $\mu$ and $\gamma$}%
\label{figure:figure-5}%
\end{figure}

For the bike sharing system, we take that $C=30$, $K=50$, $\omega=1$ and
$\mu=4$. The right one of Figure 5 indicates how the probability $p_{0}$
depends on $\lambda\in\left(  5,15\right)  $ when $\gamma=0.05,0.5$ and $1$,
respectively. It is seen that $p_{0}$ increases as $\lambda$ increases or as
$\gamma$\ decreases.

\subsection{Analysis of $p_{K}$}

Different from $p_{0}$ given in Subsection 7.1, $p_{K}$ is a probability that
the bikes are full in a tagged station, thus $p_{K}$ is also the probability
that the bike-riding customer can not return his bike at the tagged station.
To design a better bike sharing system, we hope that the value of $p_{K}$ is
as small as possible through taking a suitable parameters: $C,K,\lambda
,\mu,\gamma$ and $\omega$.

In this bike sharing system, we take that $C=30$, $K=50$, $\omega=1$ and
$\gamma=0.25$. The left one of Figure 6 shows how the probability $p_{K}$
depends on $\lambda\in\left(  10,30\right)  $ when $\mu=4,8$ and $12$,
respectively. It is seen that $p_{K}$ decreases either as $\lambda$ increases
or as $\mu$ decreases. Note that what $\lambda$ increases speeds up the rental
rate of bikes at the tagged station, while what $\mu$ decreases reduces the
return rate of bikes at the tagged station.

\begin{figure}[ptb]
\centering              \includegraphics[width=7cm]{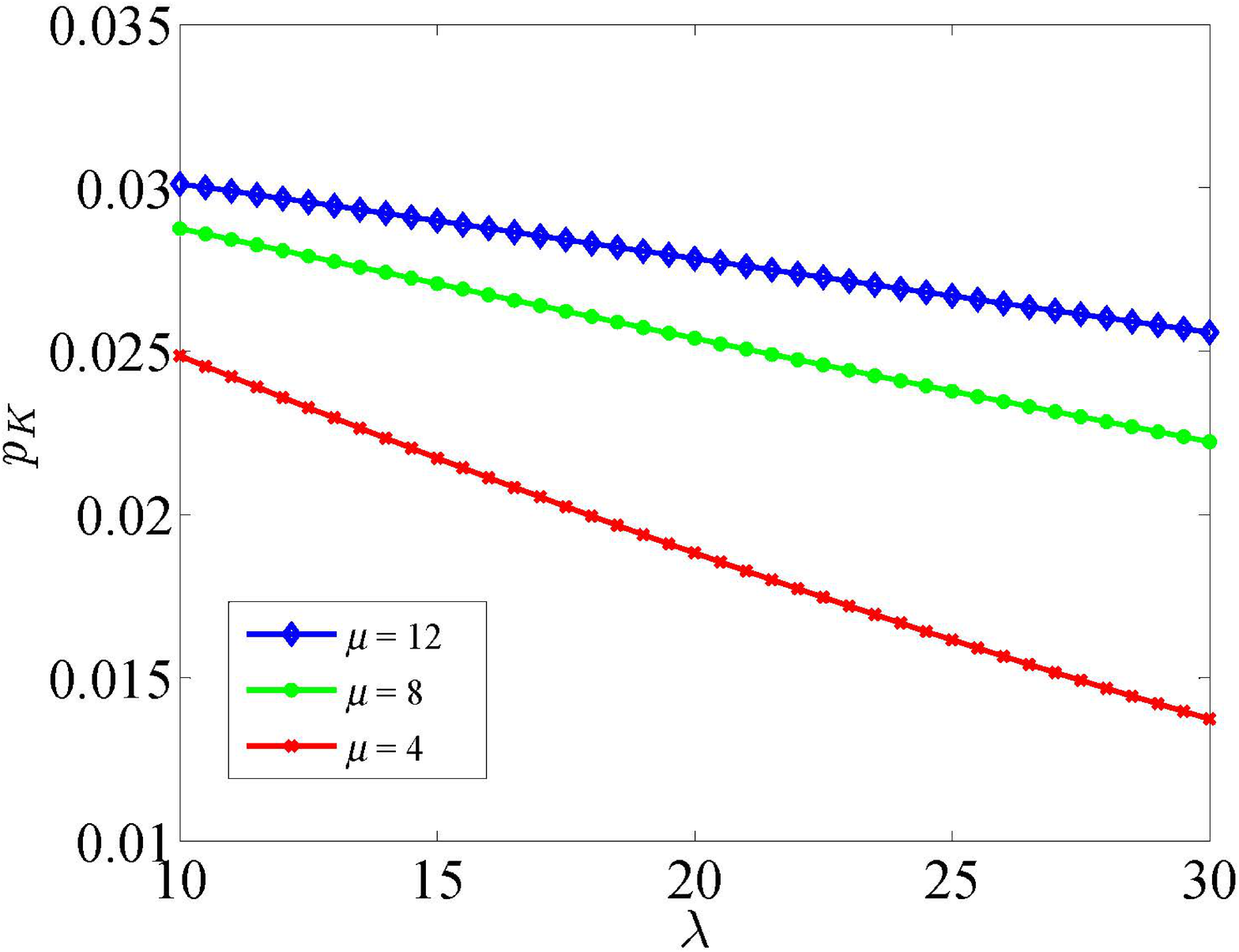}
\includegraphics[width=7cm]{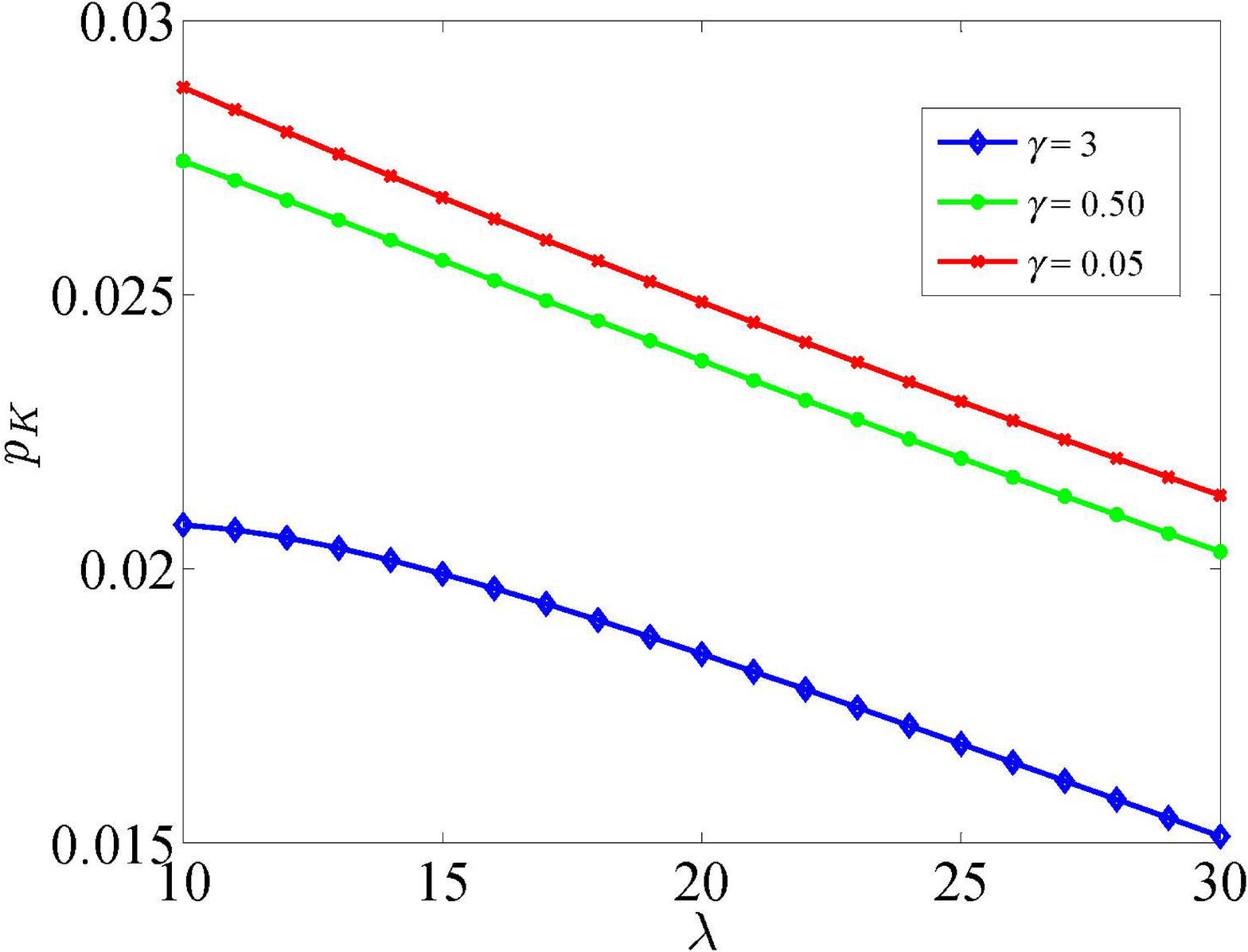}  \newline \caption{$p_{K}$ vs.
$\lambda$, $\mu$ and $\gamma$}%
\label{figure:figure-6}%
\end{figure}

For the bike sharing system, we take that $C=30$, $K=50$, $\omega=1$ and
$\mu=7$. The right one of Figure 6 indicates how the probability $p_{K}$
depends on $\lambda\in\left(  10,30\right)  $ when $\gamma=0.05,0.5$ and $3$,
respectively. It is seen that $p_{K}$ decreases as $\lambda$ increases or as
$\gamma$\ increases.

\subsection{Analysis of $p_{0}+p_{K}$}

Based on the above two analysis for $p_{0}$ and $p_{K}$, we further hope that
the value of $p_{0}+p_{K}$ can be as small as possible through taking a
suitable parameters: $C,K,\lambda,\mu,\gamma$ and $\omega$.

In this bike sharing system, we take that $C=30$, $K=50$, $\omega=1$ and
$\gamma=0.25$. The left one of Figure 7 shows how the probability $p_{0}%
+p_{K}$ depends on $\lambda\in\left(  10,30\right)  $ when $\mu=6,8$ and $10$,
respectively. It is seen that $p_{0}+p_{K}$ decreases either as $\lambda$
increases or as $\mu$ decreases. Comparing Figure 7 with Figures 5 and 6, it
is seen that $p_{K}$ has a bigger influence on the probability $p_{0}+p_{K}$
than $p_{0}$.

\begin{figure}[ptb]
\centering      \includegraphics[width=7cm]{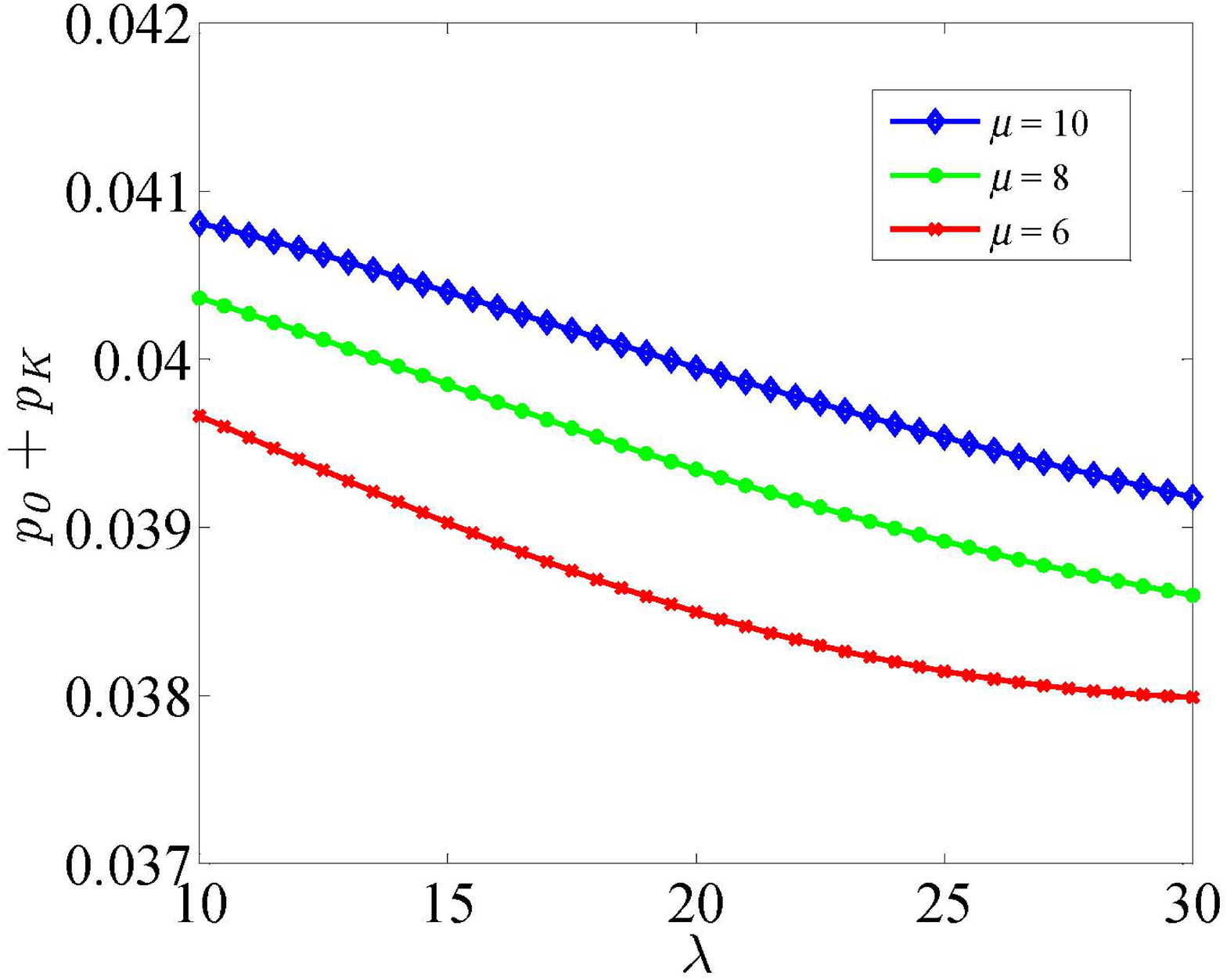}
\includegraphics[width=7cm]{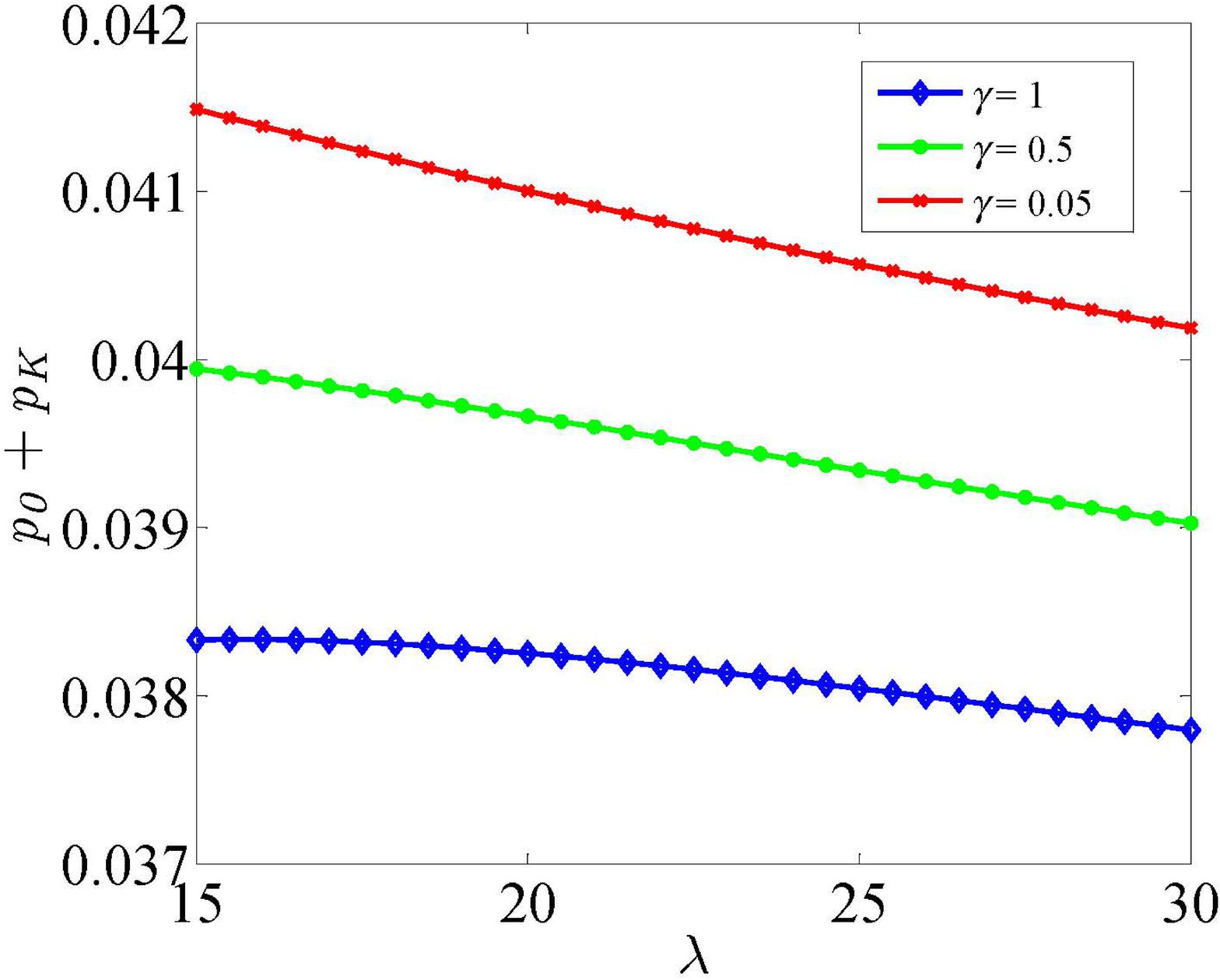}  \newline \caption{$p_{0}+p_{K}$ vs.
$\lambda$, $\mu$ and $\gamma$}%
\label{figure:figure-7}%
\end{figure}

For the bike sharing system, we take that $C=30$, $K=50$, $\omega=1$ and
$\mu=12$. The right one of Figure 7 indicates how the probability $p_{0}%
+p_{K}$ depends on $\lambda\in\left(  15,30\right)  $ when $\gamma=0.05,0.5$
and $1$, respectively. It is seen that $p_{0}+p_{K}$ decreases as $\lambda$
increases or as $\gamma$\ increases.

\subsection{Analysis of $E\left[  Q\right]  $}

From $E\left[  Q\right]  =\sum_{k=1}^{K}kp_{k}$, it is seen that $E\left[
Q\right]  $ is the stationary expected number of bikes parked at the tagged
station. Obviously, a customer who is renting a bike likes a bigger $E\left[
Q\right]  $, while a customer who is returning a bike likes a smaller
$E\left[  Q\right]  $. In addition, $E\left[  Q\right]  $ can also be used to
express the profit of the tagged station as follows:%
\[
R=-cE\left[  Q\right]  +\psi\left\{  C-E\left[  Q\right]  \right\}  ,
\]
where $c$ is the cost price per bike and per time unit when a bike is parked
in the tagged station, and $\psi$ is the benefit price per bike and per time
unit when a bike is rented from the tagged station.

In this bike sharing system, we take that $C=30$, $K=50$, $\omega=1$ and
$\gamma=0.25$. The left of Figure 8 shows how the stationary mean $E\left[
Q\right]  $ depends on $\lambda\in\left(  10,30\right)  $ when $\mu=2$, $5$
and $8$, respectively. It is seen that $E\left[  Q\right]  $ decreases either
as $\lambda$ increases or as $\mu$ decreases.

\begin{figure}[ptb]
\centering                  \includegraphics[width=7cm]{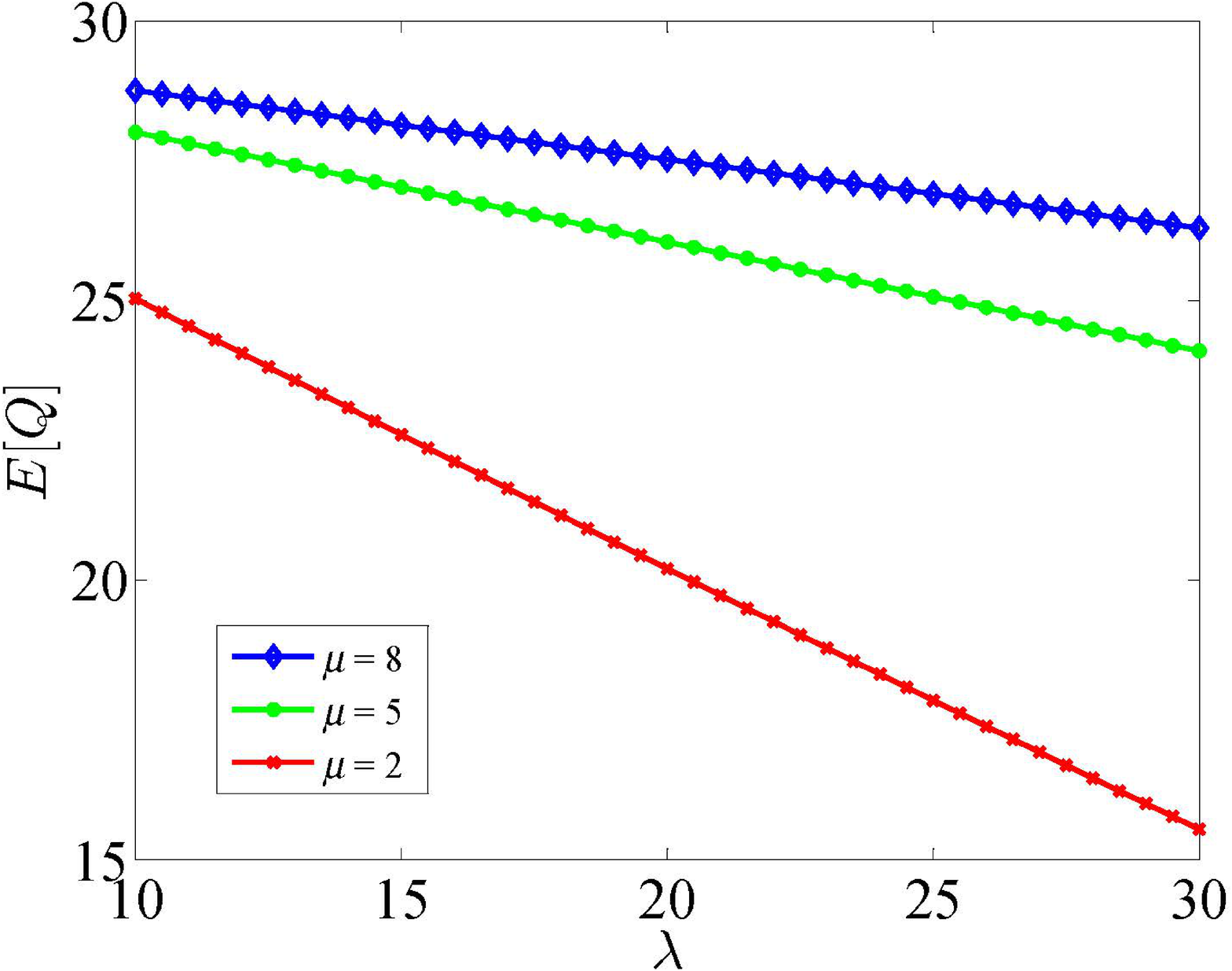}
\includegraphics[width=7cm]{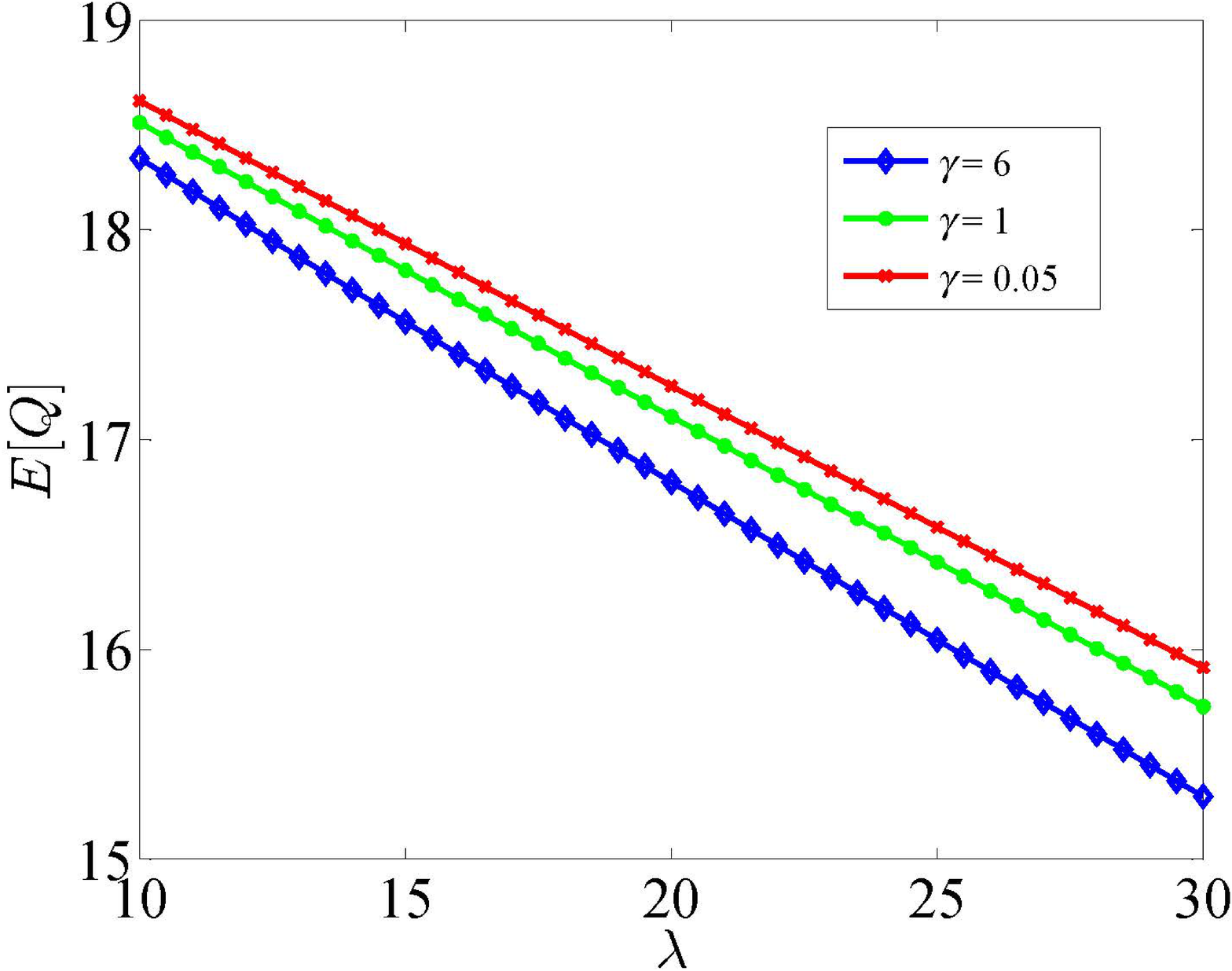}  \newline \caption{$E\left[  Q\right]
$ vs. $\lambda$, $\mu$ and $\gamma$}%
\label{figure:figure-8}%
\end{figure}

For the bike sharing system, we take that $C=20$, $K=50$, $\omega=1$ and
$\mu=7$. The right of Figure 8 indicates how the stationary mean $E\left[
Q\right]  $ depends on $\lambda\in\left(  10,30\right)  $ when $\gamma=0.05$,
$0.1$ and $6$, respectively. It is seen that $E\left[  Q\right]  $ decreases
as $\lambda$ increases or as $\gamma$ increases.

\subsection{Parameter optimization}

We provide a simple discussion for how to optimize some key parameters of the
bike sharing system through numerical experiments. Note that $\lambda$,
$\gamma$ and $\omega$\ are the arrival and walk information of any customer
respectively, thus our parameter optimization will not consider them. In this
case, our decision variables in the following optimal problems will mainly
focus on the three parameters: $C$, $K$ and $\mu$.

\textbf{(a) Optimization based on the probabilities }$p_{0}$\textbf{ and
}$p_{K}$

Since our purpose is to minimize either $p_{0}$, $p_{K}$ or $p_{0}+p_{K}$, we
may choose a weighted method in which $\beta_{1},\beta_{2}$ and $\beta_{3}$
are the weighted coefficients with $\beta_{1},\beta_{2},\beta_{3}\geq0$ and
$\beta_{1}+\beta_{2}+\beta_{3}=1$. In this case, our optimal problem is given
by%
\begin{align*}
&  \min\left\{  \beta_{1}p_{0}+\beta_{2}p_{K}+\beta_{3}\left(  p_{0}%
+p_{K}\right)  \right\} \\
\text{s.t. }  &  0<\gamma<\mu,\\
&  1\leq C\leq K.
\end{align*}
For example, when $\beta_{2}=0$ and $\beta_{3}=0$, $\min\left\{  \beta
_{1}p_{0}+\beta_{2}p_{K}+\beta_{3}\left(  p_{0}+p_{K}\right)  \right\}
=\min\left\{  p_{0}\right\}  $; when $\beta_{1}=0$ and $\beta_{3}=0$,
$\min\left\{  \beta_{1}p_{0}+\beta_{2}p_{K}+\beta_{3}\left(  p_{0}%
+p_{K}\right)  \right\}  =\min\left\{  p_{K}\right\}  $; when $\beta_{1}=0$
and $\beta_{2}=0$, $\min\left\{  \beta_{1}p_{0}+\beta_{2}p_{K}+\beta
_{3}\left(  p_{0}+p_{K}\right)  \right\}  =\min\left\{  p_{0}+p_{K}\right\}
$. Therefore, our above optimal problem is a more general tradeoff among three
key factors: $p_{0}$, $p_{K}$ and $p_{0}+p_{K}$.

\textbf{(b) Optimization based on the profit }$R$

Now, our optimal purpose is to maximize the profit of the tagged station as
follows:%
\begin{align*}
&  \max\left\{  -cE\left[  Q\right]  +\psi\left\{  C-E\left[  Q\right]
\right\}  \right\} \\
\text{s.t. }  &  0<\gamma<\mu,\\
&  1\leq C\leq K.
\end{align*}

\section{Concluding Remarks}

In this paper, we apply the mean-field theory to studying a large-scale bike
sharing system, where the mean-field computation can partly overcome the
difficulty of state space explosion in more complicated bike sharing systems.
We first use an $N$-dimensional Markov process to express the states of the
bike sharing system, and construct an empirical measure Markov process of the
$N$-dimensional Markov process. Then we set up the system of mean-field
equations by means of a virtual time-inhomogeneous $M(t)/M(t)/1/K$ queue whose
arrival and service rates are determined through some mean-field computation.
Furthermore, we employ the martingale limit to investigate the limiting
behavior of the empirical measure process, and prove that the fixed point is
unique. This illustrates the asymptotic independence of the bike sharing
system. Based on this, we can compute the fixed point through a nonlinear
birth-death process, and provide some effective algorithms for computing the
steady-state probability of the problematic stations. Finally, we use some
numerical examples to give valuable observation on how the steady-state
probability of the problematic stations depends on some crucial parameters of
the bike sharing system.

This paper provides a complete picture on how to use the mean-field theory,
the time-inhomogeneous queues, the martingale limits and the nonlinear Markov
processes to analyze performance measures of the large-scale bike sharing
systems. This picture is described as the following four key steps: (1)
Setting up system of mean-field equations, (2) proofs of the mean-field limit,
(3) uniqueness and computation of the fixed point, and (4) performance
analysis of the bike sharing system. Therefore, the methodology and results of
this paper give new highlight on understanding influence of system key
parameters on performance measures of the bike sharing systems. Along such a
line, there are a number of interesting directions for potential future
research, for example:

\begin{itemize}
\item Analyzing impact of the intelligent information technologies on
operations management of the bike sharing systems;

\item discussing the bike sharing systems with non-exponential distributions
and non-Poisson point processes, and develop some more general mean-field models;

\item studying the periodical or time-inhomogeneous bike sharing systems; and

\item modeling a bike sharing system with multiple clusters, where the
unbalanced bikes can be redistributed among the stations or clusters by means
of optimal scheduling of trucks.
\end{itemize}

\section*{Acknowledgements}

The authors thank the Area Editor and the two reviewers for many valuable
comments to sufficiently improve the presentation of this paper, and
appreciate Professor Yunan Liu at North Carolina State University for many
constructive discussions in the study of bike sharing systems. At the same
time, the first author acknowledges that this research is partly supported by
the National Natural Science Foundation of China under grant No. 71271187, No.
71471160 and No. 71671158, and the Fostering Plan of Innovation Team and
Leading Talent in Hebei Universities under grant No. LJRC027.

\end{document}